\documentclass[9pt]{article}
\usepackage{epsf}
\usepackage{epsfig}
\usepackage{url}

\def\RR{\rm \hbox{I\kern-.2em\hbox{R}}}
\def\NN{\rm \hbox{I\kern-.2em\hbox{N}}}
\def\ZZ{\rm {{\rm Z}\kern-.28em{\rm Z}}}
\def\CC{\rm \hbox{C\kern -.5em {\raise .32ex \hbox{$\scriptscriptstyle
|$}}\kern-.22em{\raise .6ex \hbox{$\scriptscriptstyle |$}}\kern .4em}}
\def\vp{\varphi}
\def\<{\langle}
\def\>{\rangle}

\def\e{\varepsilon}
\def\sm{\setminus}
\def\nl{\newline}
\def\o{\overline}

\def\bq{{\bf q}}
\def\cT{{\cal T}}
\def\cA{{\cal A}}
\def\cI{{\cal I}}

\def\cB{{\cal B}}
\def\cR{{\cal R}}
\def\cS{{\cal S}}
\def\cD{{\cal D}}
\def\cL{{\cal L}}
\def\cC{{\cal C}}
\def\cP{{\cal P}}

\def\cO{{\cal O}}
\def\Chi{\raise .3ex
\hbox{\large $\chi$}} \def\vp{\varphi}
\def\lsima{\hbox{\kern -.6em\raisebox{-1ex}{$~\stackrel{\textstyle<}{\sim}~$}}\kern -.4em}
\def\lsim{\hbox{\kern -.2em\raisebox{-1ex}{$~\stackrel{\textstyle<}{\sim}~$}}\kern -.2em}
\def\({\Bigl (}
\def\){\Bigr )}
\def\({\Bigl (}
\def\){\Bigr )}

\newcommand{\be}{\begin{equation}}
\newcommand{\ee}{\end{equation}}
\newcommand{\bea}{$$ \begin{array}{lll}}
\newcommand{\eea}{\end{array} $$}
\newcommand{\bi}{\begin{itemize}}
\newcommand{\ei}{\end{itemize}}
\newcommand{\iref}[1]{{\rm (\ref{#1})}}
\newtheorem{theorem}{Theorem}[section]
\newtheorem{remark}[theorem]{Remark}
\newtheorem{lemma}[theorem]{Lemma}

\def\cO{\mathcal O}
\def\R{\mathbb R}

\def\b{\mathbf}

\newcommand\trans{\mathrm T}

\def\gsim{\hbox{\kern -.2em\raisebox{-1ex}{$~\stackrel{\textstyle>}{\sim}~$}}\kern -.2em}

\usepackage{amsmath}
\usepackage{amsfonts}
\usepackage{amssymb}

\usepackage{psfrag}
\usepackage{graphicx}

\DeclareMathOperator{\argmin}{Argmin}

\input xy
\xyoption{all}

\begin{document}
\title{\bf Adaptive multiresolution analysis\\
based on anisotropic triangulations} 
\author{Albert Cohen,  Nira Dyn, Fr\'ed\'eric Hecht\\
and Jean-Marie Mirebeau \thanks{%
This research was supported by the P2R French-Israeli program
``Nonlinear multiscale methods - applications to image
 and terrain data''
}}
\maketitle
\date{}
\begin{abstract}
A simple greedy refinement procedure for the generation of data-adapted
triangulations is proposed and studied. Given a function $f$
of two variables, the algorithm produces a hierarchy of triangulations
$(\cD_j)_{j\geq 0}$ and piecewise polynomial approximations 
of $f$ on these triangulations. The refinement procedure consists in bisecting
a triangle $T$ in a direction which is chosen so as to minimize the local approximation error
in some prescribed norm between $f$ and its piecewise polynomial approximation
after $T$ is bisected. The hierarchical structure allows us to derive
various approximation tools such as multiresolution analysis, wavelet bases,
adaptive triangulations based either on greedy or optimal
CART trees, as well as a simple encoding of the corresponding triangulations.
We give a general proof of convergence in the $L^p$ norm
of all these approximations.
Numerical tests performed in the case of piecewise linear
approximation of functions with analytic expressions or of numerical images 
illustrate the fact that the refinement procedure generates triangles with an optimal aspect ratio (which
is dictated by the local Hessian of $f$ in case of $C^2$ functions).
\end{abstract}

\section{Introduction}
\label{intro}

Approximation by piecewise polynomial functions
is a standard procedure which occurs in various applications. 
In some of them such as 
terrain data simplification or image compression,
the function to be approximated might be fully known,
while it might be only partially known or fully
unknown in other applications such as denoising, statistical learning or 
in the finite element discretization of PDE's.

In all these applications, one usually makes the 
distinction between {\it uniform} and {\it adaptive} approximation.
In the uniform case, the domain of interest is decomposed into
a partition where all elements have comparable shape and size,
while these attributes are allowed to vary strongly in the adaptive case.
In the context of adaptive triangulations, another important distinction
is between {\it isotropic} and {\it anisotropic} triangulations.
In the first case the triangles
satisfy a condition which guarantees that 
they do not differ too much from equilateral triangles. This can
either be stated in terms of a minimal value $\theta_0>0$ for
every angle, or by a uniform bound on the
aspect ratio
$$
\rho_T:=\frac {h_T}{r_T}
$$
of each triangle $T$ where $h_T$ and $r_T$ respectively denote
the diameter of $T$ and of its largest inscribed disc.
In the second case, which is in the scope of 
the present paper, the aspect ratio is allowed
to be arbitrarily large, i.e. long and thin triangles are allowed.
In summary, adaptive and anisotropic triangulations mean that
we do not fix any constraint on the size and shape of the triangles.

Given a function $f$ and a norm $\|\cdot\|_X$ of interest, 
we can formulate the problem
of finding the {\it optimal triangulation} for $f$ in the $X$-norm 
in two related forms:
\begin{itemize}
\item
For a given $N$ find a triangulation $\cT_N$ with $N$ triangles
and a piecewise polynomial function $f_N$ (of some fixed degree) on $\cT_N$
such that $\|f-f_N\|_X$ is minimized.
\item
For a given $\e >0$ find a triangulation $\cT_N$ with minimal
number of triangles $N$ and a piecewise polynomial function $f_N$ 
such that $\|f-f_N\|_X\leq \e$.
\end{itemize}
In this paper $X$ will be the $L^p$ norm for some arbitrary $1\leq p\leq \infty$.
The exact solution to such problems is usually out of reach both
analytically and algorithmically: even when restricting the search 
of the vertices to a finite grid, the number of possible triangulations
has combinatorial complexity and an exhaustive search is therefore
prohibitive. 

Concrete mesh generation algorithms have been
developed in order to generate in reasonable time triangulations
which are ``close'' to the above described optimal trade-off 
between error and complexity. They are typically governed
by two intuitively desirable features:
\begin{enumerate}
\item
The triangulation should {\it equidistribute} the local approximation error between 
each triangle. This rationale is typically used in local mesh refinement
algorithms for numerical PDE's
\cite{Ve}: a triangle is refined when the local approximation 
error (estimated by an a-posteriori error indicator) is large.
\item
In the case of anisotropic meshes, the local aspect ratio
should in addition be optimally adapted to 
the approximated function $f$. In the case of piecewise linear
approximation, this is achieved by imposing
that the triangles are isotropic with respect to a distorted
metric induced by the Hessian $d^2f$. 
We refer in particular to \cite{BFGLS}
where this task is executed using 
Delaunay mesh generation techniques.
\end{enumerate}

While these last algorithms fastly produce 
anisotropic meshes which are naturally
adapted to the approximated function,
they suffer from two intrinsic limitations:
\begin{enumerate}
\item
They are based on the evaluation of the Hessian $d^2f$,
and therefore do not in principle apply 
to arbitrary functions $f\in L^p(\Omega)$ for $1\leq p\leq \infty$ 
or to noisy data. 
\item
They are non-hierarchical: for $N>M$, the triangulation
$\cT_N$ is not a refinement of $\cT_M$. 
\end{enumerate}

One way to circumvent the first limitation is
to regularize the function $f$, either by projection onto a 
finite element space or by convolution by a mollifier. However
this raises the additional problem of appropriately tuning
the amount of smoothing, in particular depending on the
noise level in $f$.

The need for hierarchical triangulations is critical
in the construction of wavelet bases, which play
an important role in applications
to image and terrain data processing,
in particular data compression \cite{CDDD}.
In such applications, the multilevel structure
is also of key use for the fast encoding
of the information.
Hierarchy is also useful in the 
design of optimally converging adaptive
methods for PDE's \cite{Dor,MNS,BDD,St}.
However, all these developments are so 
far mostly restricted to isotropic refinement methods.
Let us mention that hierarchical and anisotropic
triangulations have been investigated in \cite{KP},
yet in this work the triangulations are {\it fixed in advance}
and therefore generally not adapted to the approximated function.
\nl
\nl
{\it A natural objective is therefore to design
adaptive algorithmic techniques that combine
hierarchy and anisotropy, and that apply
to any function $f\in L^p(\Omega)$, without any need for regularization.}
\nl
\nl
In this paper we propose and study a
simple {\it greedy refinement procedure}
that achieves this goal: starting from an initial 
triangulation $\cD_{0}$, the procedure
bisects every triangle from one of its vertices to the mid-point of 
the opposite segment. The choice of the vertex
is typically the one which minimizes the new approximation
error after bisection among the three options.

Surprisingly, it turns out that - in the case of piecewise
linear approximation - this elementary
strategy tends to generate anisotropic triangles 
with optimal aspect ratio. This fact is rigorously
proved in \cite{CM} which establishes optimal 
error estimates for the approximation
of {\it smooth and convex functions $f\in C^2$},
by adaptive triangulations $\cT_N$ with $N$ triangles.
These triangulations are obtained
by consecutively applying the refinement procedure 
to the triangle of maximal error. The estimates in \cite{CM}
are of the form
\be
\|f-f_N\|_{L^p}  \leq CN^{-1}\|\sqrt{|{\rm det}(d^2f)|}\|_{L^\tau},\;\; \frac 1 \tau=\frac 1 p + 1.
\label{aniser}
\ee
and were already established in \cite{CSX,BBLS} for functions which
are not necessarily convex, however based on triangulations which are non-hierarchical and based on the
evaluation of $d^2f$. Note that \iref{aniser} improves on the estimate
\be
\|f-f_N\|_{L^p}  \leq CN^{-1}\|d^2f\|_{L^\tau},\;\; \frac 1 \tau=\frac 1 p + 1.
\label{isoer}
\ee
which can be established (see \S 2 in \cite{CM}) for adaptive triangulations with 
isotropic triangles, and which itselfs improves on the
classical estimate
\be
\|f-f_N\|_{L^p}  \leq CN^{-1}\|d^2f\|_{L^p},
\label{unier}
\ee
which is known to hold for uniform triangulations.

The main objective of the present paper is to introduce
the refinement procedure as well as several approximation
methods based on it, and
to study their convergence for {\it an arbitrary function $f\in L^p$}.
In \S 2, we introduce notation that serves to
describe the refinement procedure and define
the anisotropic hierarchy of triangulations $(\cD_j)_{j\geq 0}$. 
We show how this general framework can be used 
to derive adaptive approximations of $f$ either by 
triangulations based on greedy or optimal trees,
or by wavelet thresholding.
In \S 3, we show that as defined, the approximations
produced by the refinement
procedure may fail to converge for 
certain $f \in L^p$ and show how to modify 
the procedure so that convergence
holds for any arbitrary $f\in L^p$.
We finally present  in \S 4 some numerical tests which illustrate 
the optimal mesh adaptation, in the case of piecewise linear
elements, when the refinement procedure is applied 
either to synthetic functions
or to numerical images. 

\section{An adaptive and anisotropic multiresolution framework}

\subsection{The refinement procedure}

Our refinement procedure is based on
a local approximation operator
$\cA_T$ 
acting from $L^p(T)$ onto $\Pi_m$ - the space of polynomials of
total degree less or equal to $m$. Here, the parameters $m\geq 0$ and $1\leq p\leq \infty$
are arbitrary but fixed. For a generic triangle
$T$, we denote by $(a,b,c)$
its edge vectors oriented in clockwise 
or anticlockwise direction so that
$$
a+b+c=0.
$$
We define the
local $L^p$ approximation error 
$$
e_T(f)_p:=\|f-\cA_Tf\|_{L^p(T)}.
$$
The most natural choice for $\cA_T$ is the operator $\cB_T$ of best $L^p(T)$ approximation
which is defined by
$$
\|f-\cB_Tf\|_{L^p(T)}=\min_{\pi\in\Pi_m}\|f-\pi\|_{L^p(T)}.
$$
However this operator is non-linear and not easy to compute when $p\neq 2$. 
In practice, one prefers to use an operator 
which is easier to compute, yet nearly optimal in the sense
that
\be
\|f-\cA_Tf\|_{L^p(T)} \leq C \inf_{\pi\in\Pi_m}\|f-\pi\|_{L^p(T)}.
\label{unileb}
\ee
with $C$ a Lebesgue constant independent of $f$ and $T$. 
Two particularly simple admissible choices of approximation operators are the following:
\begin{enumerate}
\item
$\cA_T=P_T$, the $L^2(T)$-orthogonal projection onto $\Pi_m$, defined
by $P_Tf\in \Pi_m$ such that ${\int_T(f-P_T f)\pi}=0$ for all $\pi\in\Pi_m$. This operator 
has finite Lebesgue constant for all $p$, with $C=1$ when $p=2$
and $C>1$ otherwise.
\item
$\cA_T=I_T$, the local interpolation operator which is defined by  $I_Tf\in\Pi_m$ such that
$I_Tf(\gamma)=f(\gamma)$ for all $\gamma\in \Sigma:=\{ \sum \frac {k_i}m v_i\; ;Ê\; k_i\in \NN,\; 
\sum k_i=m\}$
where $\{v_1,v_2,v_3\}$ are the vertices of $T$ (in the case $m=0$ we can take for $\Sigma$
the barycenter of $T$). This operator is only defined on
continuous functions and has Lebesgue constant $C>1$ in the $L^\infty$ norm. 
\end{enumerate}
All our results are simultaneously valid when $\cA_T$ is either $P_T$
or $I_T$ (in the case where $p=\infty$), or any linear operator that fulfills 
the continuity assumption \iref{unileb}.

Given a target function, our refinement procedure defines by induction a hierarchy of nested
triangulations  $(\cD_j)_{j\geq 0}$ with $\#(\cD_j)=2^j \#(\cD_0)$.
The procedure starts from the coarse triangulation $\cD_0$ of $\Omega$,
which is fixed independently of $f$. When 
$\Omega=[0,1]^2$ we may split it into
two symmetric triangles so that $\#(\cD_0)=2$. 
For every $T\in\cD_j$, we split $T$ into 
two sub-triangles of equal area by bisection
from one of its three vertices towards the
mid-point of the opposite edge $e\in \{a,b,c\}$. 
We denote by 
$T_e^1$ and $T_e^2$
the two resulting triangles. The choice of $e\in \{a,b,c\}$
is made according to a {\it refinement rule}
that selects this edge depending on the properties of $f$.
We denote by $\cR$ this refinement rule, which can therefore
be viewed as a mapping
$$
\cR : (f,T) \mapsto e.
$$
We thus obtain two {\it children} of $T$ corresponding to the choice $e$.
$\cD_{j+1}$ is the triangulation consisting of all such pairs 
corresponding to all $T\in \cD_j$.

In this paper, we consider refinement rules
where the selected edge $e$ minimizes a {\it decision function} $e\mapsto d_T(e,f)$ 
among $\{a,b,c\}$. We refer to such rules
as {\it greedy refinement rules}. A more elaborate
type of refinement rule is also considered in \S 3.3.

The role of the decision function is to drive the generation of 
anisotropic triangles according to the local properties of
$f$, in contrast to simpler procedures
such as {\it newest vertex bisection} (i.e. split $T$ from the most recently created vertex)
which is independent of $f$ and generates triangulations with isotropic
shape constraint. 

Therefore, the choice of $d_T(e,f)$ is critical in
order to obtain triangles with an optimal aspect ratio.
The most natural choice corresponds to the optimal split
$$
d_T(e,f)=e_{T_e^1}(f)_p^p+e_{T_e^2}(f)_p^p,
$$
i.e. choose the edge that minimizes the resulting $L^p$ error after bisection. 
It is proved in \cite{CM} in the case of piecewise linear
approximation, that when $f$ is
a $C^2$ function which is strictly convex or concave
the refinement rule based on the decision function 
\be
d_T(e,f)=\|f-I_{T_e^1}f\|_{L^1(T_e^1)}+\|f-I_{T_e^2}f\|_{L^1(T_e^2)}.
\label{optil1}
\ee
generates triangles which tend to have have an optimal aspect ratio,
locally adapted to the Hessian $d^2f$. This aspect ratio
is independent of the $L^p$ norm in which one wants to minimize the error
between $f$ and its piecewise affine approximation.

\begin{remark}
If the minimizer $e$ is not unique, we may choose it among 
the multiple minimizers either randomly 
or according to some prescribed ordering of the edges (for example
the largest coordinate pair of the opposite vertex in lexicographical order).
\end{remark} 

\begin{remark}
The triangulations $\cD_j$ which are generated by
the greedy procedure are in general non-conforming, i.e.
exhibit hanging nodes. This is not problematic in the present setting since
we consider approximation in the $L^p$ norm which does
not require global continuity of the piecewise polynomial functions.
\end{remark}

The refinement rule $\cR$ defines a multiresolution
framework. For a given $f\in L^p(\Omega)$ and any 
triangle $T$ we denote by
$$
\cC(T):=\{T_1,T_2\},
$$
the {\it children} of $T$ which are 
the two triangles obtained by splitting $T$
based on the prescribed decision function $d_T(e,f)$. We also
say that $T$ is the parent of $T_1$ and $T_2$ and
write
$$
T=\cP(T_1) = \cP(T_2).
$$
Note that 
$$
\cD_j:=\cup_{T\in\cD_{j-1}}\cC(T).
$$
We also define
$$
\cD:=\cup_{j\geq 0}\cD_j,
$$
which has the structure of {\it an infinite binary tree}.
Note that $\cD_j$ depends on $f$ (except for $j=0$)
and on the refinement rule $\cR$,
and thus $\cD$ also depends on $f$ and $\cR$:
$$
\cD_j=\cD_j(f,\cR)\;\; {\rm and}\;\; \cD=\cD(f,\cR).
$$
For notational simplicity, we sometimes omit the dependence
in $f$ and $\cR$ when there is
no possible ambiguity.

\subsection{Adaptive tree-based triangulations}

A first application of the multiresolution framework
is the design of adaptive anisotropic triangulations $\cT_N$
for piecewise polynomial approximation, by a {\it greedy 
tree algorithm}. For any finite sub-tree $\cS\subset \cD$, we denote
by 
$$
\cL(\cS):=\{T\in \cS\;\; {\rm s.t.}\;\; \cC(T)\notin\cS\}
$$
its {\it leaves} which form a partition of $\Omega$. We also denote by
$$
\cI(\cS):=\cS\sm\cL(\cS),
$$
its {\it inner nodes}. Note that {\it any} finite partition 
of $\Omega$ by elements of $\cD$ 
is the set of leaves of a finite sub-tree. One easily checks that 
$$
\#(\cS)=2\#(\cL(\cS))-N_0.
$$
For each $N$, the greedy tree algorithm defines a finite sub-tree $\cS_N$
of $\cD$ which grows from $\cS_{N_0}:=\cD_0=\cT_{N_0}$, by adding
to $\cS_{N-1}$ the two children of the triangle $T^*_{N-1}$
{\it which maximizes the local $L^p$-error
$e_T(f)_p$ over all triangles in $\cT_{N-1}$}.

The adaptive partition $\cT_N$ associated
with the greedy algorithm is defined by
$$
\cT_N:=\cL(\cS_N).
$$
Similarly to $\cD$, the triangulation $\cT_N$ depends on $f$ and on
the refinement rule $\cR$, but also on $p$
and on the choice of the approximation operator $\cA_T$.
We denote by $f_N$ the 
piecewise polynomial approximation to $f$ which is defined as $\cA_T f$ on each $T\in \cT_N$. The global $L^p$ approximation
error is thus given by
$$
\|f-f_N\|_{L^p}=\|(e_T(f)_p)\|_{\ell^p(\cT_N)}.
$$
Stopping criterions for the algorithm
can be defined in various ways:
\begin{itemize}
\item
Number of triangles: stop once a prescribed
$N$ is attained.
\item
Local error: stop once $e_T(f)_p\leq \e$ for all $T\in\cT_N$, for some prescribed $\e>0$.
\item
Global error: stop once $\|f-f_N\|_{L^p}\leq \e$ for some prescribed $\e>0$.
\end{itemize}

\begin{remark}
The role of the
triangle selection based on the largest $e_T(f)_p$ is to
equidistribute the local $L^p$ error, a feature which is desirable
when we want to approximate $f$ in $L^p(\Omega)$
with the smallest number of triangles.
However, it should be well understood that the
refinement rule may still be chosen 
based on a decision function defined
by approximation errors in norms that differ from $L^p$.
In particular, as explained earlier, the decision 
function \iref{optil1} generates triangles 
which tend to have have an optimal aspect ratio,
locally adapted to the Hessian $d^2f$ when $f$ is strictly convex
or concave, and this aspect ratio
is independent of the $L^p$ norm in which one wants to minimize the error
between $f$ and its piecewise affine approximation.
\end{remark}

The greedy algorithm is one particular
way of deriving an adaptive triangulation
for $f$ within the multiresolution framework
defined by the infinite tree $\cD$. An interesting
alternative is 
to build
adaptive triangulations within $\cD$ which
offer an {\it optimal trade-off between 
error and complexity}. This can be done
when $1\leq p<\infty$ by 
solving the minimization problem
\be
\min_{\cS} \Bigl\{ \sum_{T\in\cL(\cS)} e_T(f)_p^p + \lambda \#(\cS) \Bigl\}
\label{mincart}
\ee
among all finite trees, for some fixed $\lambda>0$.
In this approach, we do not directly control the
number of triangles which depends on the penalty
parameter $\lambda$. However, it is immediate
to see that if $N=N(\lambda)$ is the cardinality of 
$\cT^*_N=\cL(\cS^*)$ where $\cS^*$ is the 
minimizing tree, then $\cT^*_N$
minimizes the $L^p$ approximation error
$$
\cT^*_N:=\underset {\#(\cT)\leq N}\argmin\sum_{T\in\cT}e_T(f)_p^p,
$$
where the minimum is taken among all partitions $\cT$ of $\Omega$ 
within $\cD$ of cardinality less than or equal to $N$.

Due to the additive structure of the error term, the
minimization problem \iref{mincart} can be performed
in fast computational time using an optimal pruning algorithm
of CART type, see \cite{BFOS,Do}.
In the case $p=\infty$ the associated minimization problem
\be
\min_{\cS} \Bigl\{ \sup_{T\in\cL(\cS)} e_T(f)_\infty + \lambda \#(\cS)\Bigl\},
\label{mincartinf}
\ee
can also be solved by a similar fast algorithm.
It is obvious that this method improves over
the greedy tree algorithm: 
if $N$ is the cardinality of the
triangulation resulting from the minimization in \iref{mincart}
and $f_N^*$ the corresponding piecewise polynomial 
approximation of $f$ associated with this triangulation,
we have
$$
\|f-f_N^*\|_{L^p}\leq \|f-f_N\|_{L^p},
$$
where $f_N$ is built by the greedy tree algorithm.

\subsection{Anisotropic wavelets}

The multiresolution framework allows us to 
introduce the piecewise polynomial multiresolution spaces
$$
V_j=V_j(f,\cR):=\{g \;\; {\rm s.t.}\;\; g_{|T}\in \Pi_m,\; T\in \cD_j\},
$$
which depend on $f$ and on the refinement rule $\cR$. These spaces
are nested and we denote  by
$$
V=V(f,\cR)=\cup_{j\geq 0}V_j(f,\cR),
$$
their union. For notational simplicity, we sometimes omit the dependence
in $f$ and $\cR$ when there is
no possible ambiguity.

The $V_j$ spaces may be used
to construct wavelet bases, following the approach
introduced in \cite{Alp}
and that we describe in our present setting.

The space $V_j$ is equipped with an orthonormal
{\it scaling function basis}:
$$
\vp_T^{i},\;\; \; i=1,\cdots,\frac 1 2(m+1)(m+2),\;\; T\in\cD_j, 
$$
where the $\vp_T^i$ 
for $i=1,\cdots,\frac 1 2(m+1)(m+2)$
are supported in $T$ and constitute an orthonormal basis
of $\Pi_m$ in the sense of $L^2(T)$ for each $T\in\cD$. There
are several possible choices for such a basis. In the particular case where
$m=1$, a simple one is to take
for $T$ with vertices $(v_1,v_2,v_3)$,
$$
\vp_T^i(v_i)=|T|^{-1/2} \sqrt 3\;\;{\rm and}\;\; \vp_T^i(v_j)=-|T|^{-1/2} \sqrt 3,\; j\neq i.
$$
We denote by $P_j$ the orthogonal projection onto $V_j$:
$$
P_jg:=\sum_{T\in\cD_j} \sum_{i}\<g,\vp_T^i\>\vp_T^i.
$$
We next introduce for each $T\in\cD_j$
a set of {\it wavelets}
$$
\psi_T^{i},\;\;\; i=1,\cdots,\frac 1 2(m+1)(m+2),
$$
which constitutes an orthonormal basis
of the orthogonal complement of  $\Pi_m(T)$ into $\Pi_m(T')\oplus\Pi_m(T'')$ 
with $\{T',T''\}$ the children of $T$. 
In the particular case where
$m=1$, a simple choice for
such a basis is as follows: if $(v_1,v_2,v_3)$ and $(w_1,w_2,w_3)$
denote the vertices of $T'$ and $T''$, with the
convention that $v_1=w_1$ and $v_2=w_2$ denote
the common vertices, the second one being the midpoint of the segment $(v_3,w_3)$
(i.e. $T$ has vertices $(v_3,w_3,v_1)$), then
$$
\begin{array}{ll}
& \psi_T^1:=\frac{\vp^3_{T'}-\vp^3_{T''}}{\sqrt 2},\\
& \psi_T^2:=\frac{\vp^1_{T'}-\vp^2_{T'}-\vp^1_{T''}+\vp^2_{T''}} 2
,\\
&\psi_T^3:=\frac {\vp^1_{T'}-\vp^3_{T'}+\vp^1_{T''}-\vp^3_{T''}} 2.
\end{array}
$$
where $\vp^i_{T'}$ and $\vp^i_{T''}$ are the above defined scaling functions. 

The family 
$$
\psi_T^{i},\;\;\; i=1,\cdots,\frac 1 2(m+1)(m+2),\;\; T\in\cD_j
$$
constitutes an orthonormal basis of $W_j$, the $L^2$-orthogonal complement of
$V^j$ in $V^{j+1}$.  A multiscale orthonormal basis of $V_J$ is given by
$$
\{\vp_T\}_{T\in\cD_0}\cup \{\psi_T^i\}_{\substack{i=1,\cdots,\frac 1 2(m+1)(m+2),\\T\in\cD_j, \, j=0,\cdots,J-1}}.
$$
Letting $J$ go to $+\infty$ we thus obtain that 
$$
\{\vp_T\}_{T\in\cD_0}\cup \{\psi_T^i\}_{\substack{i=1,\cdots,\frac 1 2(m+1)(m+2),\\T\in\cD_j, \, j\geq 0}}
$$
is an orthonormal basis of the space
$$
V(f,\cR)_2:=\overline{V(f,\cR)}^{L^2(\Omega)}=\overline{\cup_{j\geq 0} V_j(f,\cR)}^{L^2(\Omega)}.
$$
For the sake of notational simplicity, we rewrite this basis as 
$$
(\psi_\lambda)_{\lambda\in\Lambda},
$$
Note that $V(f,\cR)$ is not necessarily dense
in $L^2(\Omega)$ and so $V(f,\cR)_2$ is not always equal to $L^2(\Omega)$.
Therefore, the expansion of an arbitrary function $g\in L^2(\Omega)$
in the above wavelet basis does not always converge towards $g$
in $L^2(\Omega)$. The same remark holds for the $L^p$ convergence
of the wavelet expansion of an arbitrary 
function $g\in L^p(\Omega)$ (or $\cC(\Omega)$
in the case $p=\infty$): $L^p$-convergence holds when the space
$$
V(f,\cR)_p:=\overline{V(f,\cR)}^{L^p(\Omega)}
$$
coincides with $L^p(\Omega)$ (or contains $\cC(\Omega)$ in the case $p=\infty$),
since we have
$$
\Big\|f-\sum_{|\lambda|<j}d_\lambda\psi_\lambda\Big\|_{L^p}=\|f-P_jf\|_{L^p} \leq C
\inf_{g\in V_j}\|f-g\|_{L^p},
$$
with $C$ the Lebesgue constant in \iref{unileb} for the orthogonal projector.

A sufficient condition for such a property to hold is obviously that the
size of all triangles goes to $0$ as the level $j$ increases, i.e.
$$
\lim_{j\to +\infty}\sup_{T\in\cD_j} {\rm diam}(T)=0.
$$
However, this condition might not hold for the
hierarchy $(\cD_j)_{j\geq 0}$ produced by the 
refinement procedure.
On the other hand, the multiresolution 
approximation being intrinsically adapted to $f$,
a more reasonable requirement is that the
expansion of $f$ converges towards $f$ in $L^p(\Omega)$
when $f\in L^p(\Omega)$ (or $\cC(\Omega)$
in the case $p=\infty$). This is equivalent to the property 
$$
f\in V(f,\cR)_p.
$$
We may then define an adaptive approximation of $f$
by thresholding its coefficients
at some level $\e>0$:
$$
f_\e:=\sum_{|f_\lambda|\geq \e} f_\lambda\psi_\lambda,
$$
where $f_\lambda:=\<f,\psi_\lambda\>$. 
When measuring the error in the $L^p$ norm, a more natural choice is
to perform thresholding on the component of the
expansion measured in this norm, defining therefore
$$
f_\e:=\sum_{\|f_\lambda\psi_\lambda\|_{L^p}\geq \e} f_\lambda\psi_\lambda.
$$
We shall next see that the condition $f\in V(f,\cR)_p$ also ensures
the convergence of the tree-based adaptive approximations $f_N$ and $f_N^*$ towards $f$
in $L^p(\Omega)$. We shall also see that this condition may not 
hold for certain functions $f$, but that this
difficulty can be circumvented by a modification
of the refinement procedure.

\section{Convergence analysis}

\subsection{A convergence criterion}

The following result relates the convergence towards $f$ of
its approximations by projection onto the spaces $V_j(f,\cR)$, greedy and optimal tree algorithms, 
and wavelet thresholding. This result is valid for {\it any} refinement rule $\cR$.

\begin{theorem}
Let $\cR: (f,T)\mapsto e$ be an arbitrary refinement rule and let
$f\in L^p(\Omega)$. The following statements are equivalent:

{\rm (i)} $f\in V(f,\cR)_p$.

{\rm (ii)} The greedy tree approximation converges: $\lim_{N\to +\infty}\|f-f_N\|_{L^p}=0$.

{\rm (iii)} The optimal tree approximation converges:
$\lim_{N\to +\infty}\|f-f^*_N\|_{L^p}=0$.
\nl
In the case $p=2$, they are also equivalent to:

{\rm (iv)} The thresholding approximation converges:
$\lim_{\e\to 0}\|f-f_\e\|_{L^2}=0$.
\end{theorem}
\noindent
{\bf Proof:}  Clearly, (ii) implies (iii) since $\|f-f_N^*\|_{L^p}\leq \|f-f_N\|_{L^p}$.
Since the triangulation $\cD_N$ is a refinement of $\cT^*_N$,
we also find that  (iii) implies (i) as $\inf_{g\in V_N}\|f-g\|_{L^p}\leq \|f-f^*_N\|_{L^p}$.

We next show that (i) implies (ii). We first note that 
a consequence of (i) is that 
$$
\lim_{j\to+\infty}\sup_{T\in\cD_j} e_T(f)_p =0.
$$
It follows that for any $\eta>0$, there exists
$N(\eta)$ such that for $N>N(\eta)$, all triangles $T\in\cT_N$
satisfy
$$
e_T(f)_p\leq \eta.
$$
On the other hand, (i) means that for all $\e>0$, there exists
$J=J(\e)$ such that
$$
\inf_{g\in V_J}\|f-g\|_{L^p}\leq \e.
$$
For $N>N(\eta)$, we now split $\cT_N$ into $\cT_N^+\cup\cT_N^-$ where
$$
\cT_N^+:=\cT_N\cap (\cup_{j\geq J} \cD_j)\;\;{\rm and}
\;\; \cT_N^-:=\cT_N\cap (\cup_{j< J} \cD_j).
$$
We then estimate the error of the greedy algorithm by
\begin{eqnarray*}
\|f-f_N\|_{L^p}^p & = & \sum_{T\in\cT_N^+}e_T(f)_p^p+\sum_{T\in\cT_N^-}e_T(f)_p^p\\
& \leq & C^p\sum_{T\in\cT_N^+}\inf_{\pi\in \Pi_m}\|f-\pi\|_{L^p(T)}^p+\eta^p \#(\cT_N^-) \\
& \leq & C^p\inf_{g\in V_J}\|f-g\|_{L^p}^p+\eta^p \#(\cT_N^-) \\
&\leq & C^p\e^p+2^JN_0\eta^p,
\end{eqnarray*}
%
where $C$ is the stability constant of \iref{unileb}.
This implies (ii) since for any $\delta>0$,
we can first choose $\e>0$ such that $C^p\e^p<\delta/2$, and
then choose $\eta>0$ such that $2^{J(\e)}N_0\eta^p<\delta/2$.
When $p=\infty$ the estimate is modified into
$$
\|f-f_N\|_{L^\infty}\leq \max\{C\e,\eta\},
$$
which also implies (ii) by a similar reasoning.

We finally prove the equivalence
between (i) and (iv) when $p=2$. Property (i) is  equivalent to
the $L^2$ convergence of the orthogonal projection 
$P_{j}f$ to $f$ as $j\to +\infty$, or 
equivalently of the partial sum 
$$
\sum_{|\lambda|< j} f_\lambda \psi_\lambda
$$
where $|\lambda|$ stands for the scale level of the wavelet $\psi_\lambda$.
Since $(\psi_\lambda)_{\lambda\in\Lambda}$ is an orthonormal basis of $V(f,\cR)_2$,
the summability and limit of $\sum_{\lambda\in\Lambda} f_\lambda \psi_\lambda$
do not depend on the order of the terms. Therefore (i) is equivalent
to the convergence of $f_\e$ to $f$.\hfill $\diamond$

\begin{remark}
The equivalence between statements (i) and (iv) can be extended to $1<p<\infty$
by showing that $(\psi_\lambda)_{\lambda\in\Lambda}$
is an $L^p$-unconditional system. Recall that this property means
that there exists an absolute constant $C>0$ such that
for any finitely supported sequences $(c_\lambda)$ and $(d_\lambda)$ such
$|c_\lambda|\leq |d_\lambda|$ for all $\lambda$, one has
$$
\| \sum c_\lambda\psi_\lambda\|_{L^p}\leq C\|\sum d_\lambda \psi_\lambda\|_{L^p}.
$$
A consequence of this property is that if $f\in L^p$
can be expressed as the $L^p$ limit 
$$
f=\lim_{j\to +\infty} \sum_ {|\lambda|<j} f_\lambda\psi_\lambda,
$$
then any rearrangement of the series $\sum f_\lambda \psi_\lambda$
converges towards $f$ in $L^p$. This easily implies the 
equivalence between (i) and (iv). The fact that 
$(\psi_\lambda)_{\lambda\in\Lambda}$ is an unconditional system
is well known for Haar systems \cite{KS} which correspond to the case $m=0$,
and can be extended to $m>0$ in a straightforward manner.
\end{remark}

\subsection{A case of non-convergence}

We now show that if we use a greedy refinement
rule $\cR$ based on a decision function either based
on the interpolation or $L^2$ projection error
after bisection, there exists functions $f\in L^p(\Omega)$ such that
$f\notin V(f,\cR)_p$. Without loss of generality, it is enough to construct $f$
on the reference triangle $T_{\rm ref}$ of vertices $\{(0,0),(1,0),(1,1)\}$, since our construction
can be adapted to any triangle by an affine change of variables.

Consider first a decision function defined
from the interpolation error after bisection, such as \iref{optil1}, or
more generally 
$$
d_T(f,e):=\|f-I_{T_e^1}f\|_{L^p(T_e^1)}^p+\|f-I_{T_e^2}f\|_{L^p(T_e^2)}^p.
$$
Let $f$ be a continuous function which is not identically $0$
on $T_{\rm ref}$ and which vanishes
at all points $(x,y)$ such that $x=\frac k{2m}$ for $k=0,1,\cdots,2m$,
where $m$ is the degree of polynomial approximation.
For such an $f$, it is easy to see that
$I_Tf=0$ and that $I_{T'}f=0$ for all
subtriangles $T'$ obtained by one
bisection of $T_{\rm ref}$. This shows that there is 
no preferred bisection. Assuming
that we bisect from the vertex $(0,0)$ to the opposite mid-point $(1,\frac 1 2)$, we 
find that a similar situation occurs when splitting the two subtriangles.
Iterating this observation, we see that an admissible choice of bisections
leads after $j$ steps to a triangulation $\cD_j$ of $T_{\rm ref}$ consisting of the triangles $T_{j,k}$ with vertices
$\{(0,0),(1,2^{-j}k),(1,2^{-j}(k+1))\}$ with $k=0,\cdots,2^j-1$. On each of these triangles
$f$ is interpolated by the null function and therefore
by \iref{unileb} the best $L^\infty$ approximation in $V_j$ does not converge
to $f$ as $j\to +\infty$, i.e.  $f\notin V(f,\cR)_\infty$. It can also easily be
checked that $f\notin V(f,\cR)_p$.

Similar counter-examples can be constructed
when the decision function is defined
from the $L^2$ projection error, and has the form
$$
d_T(f,e):=\|f-P_{T_e^1}f\|_{L^p(T_e^1)}^p+\|f-P_{T_e^2}f\|_{L^p(T_e^2)}^p.
$$
Here, we describe such a construction in the case $m=1$. We define
$f$ on $\cR$ as a function of the first variable given by
$$
f(x,y)=u(x),\;\;{\rm if}\;\;x\in \left[0,\frac 1 2\right], \;\; f(x,y)=u\left(x-\frac 1 2\right),\;\;{\rm if}\;\;x\in \left(\frac 1 2,1\right],
$$
where $u$ is a non-trivial function in $L^2([0,\frac 1 2])$ such that
$$
\int_{0}^{\frac 1 2} u(x)dx=\int_{0}^{\frac 1 2} xu(x)dx=\int_{0}^{\frac 1 2} x^2u(x)dx=0.
$$
A possible choice is $u(x)=L_3(4x-1)=160x^3-120x^2+24x-1$ where $L_3$ is the Legendre polynomial
of degree $3$ defined on $[-1,1]$. With this choice, we have the following result.

\begin{lemma}
Let $T$ be any triangle such that its
vertices have $x$ coordinates either 
$(0,\frac 1 2,1)$ or $(0,1,1)$ or $(\frac 1 2,1,1)$.
Then $f$ is orthogonal to $\Pi_1$ in $L^2(T)$.
\end{lemma}
\noindent
{\bf Proof:} Define
$$
T_0:=\left\{(x,y)\in T,\; x\in \left[0,\frac 1 2\right]\right\}\;\;{\rm and}\;\; T_1:=\left\{(x,y)\in T,\; x\in \left(\frac 1 2,1\right]\right\}.
$$
Then, with $v(x,y)$ being either the function $1$ or $x$ or $y$, we have
$$
\begin{array}{ll}
\int_T f(x,y)v(x,y) dxdy&=\int_{T_0}u(x)v(x,y)dxdy +\int_{T_1}u(x)v(x,y)dxdy \\
& =\int_{0}^{\frac 1 2}u(x)q_0(x)dx +\int_{0}^{\frac 1 2}u(x)q_1(x)dx,
\end{array}
$$
with
$$
q_0(x):=\int_{T_{0,x}} v(x,y)dy,\;\; q_1(x):=\int_{T_{1,x}}v(x,y)dy,
$$
where $T_{i,x}=\{y\; :\; (x,y)\in T_i\}$ for $i=0,1$.
The functions $q_0$ and $q_1$ are polynomials of degree at most $2$
and we thus obtain from the properties of $u$ that  $\int_T fv=0$. \hfill $\diamond$
\nl
\nl
The above lemma shows that for any of the three possible choices
of bisection of $T_{\rm ref}$ based on the $L^2$ decision function, 
the error is left unchanged since the projection
of $f$ on all possible sub-triangle is $0$. There is 
therefore no preferred choice,
and assuming that we bisect from the vertex $(0,0)$ to the opposite mid-point $(1,\frac 1 2)$, then
we see that a similar situation occurs when splitting the two subtriangles.
The rest of the arguments showing that $f\notin V(f,\cR)_p$
are the same as in the previous counter-example.

The above two examples of non-convergence reflect the fact 
that when $f$ has some {\it oscillations}, the refinement procedure
cannot determine the most appropriate bisection. In order
to circumvent this difficulty one needs to modify
the refinement rule.

\subsection{A modified refinement rule}

Our modification consists of bisecting from the most recently generated
vertex of $T$, in case the local error is not
reduced enough by all three bisections.  More precisely, 
we modify the choice of the bisection of any $T$ as follows:

Let $e$ be the edge which minimizes the decision function $d_T(e,f)$.
If 
$$
\left(e_{T_e^1}(f)_p^p+e_{T_e^2}(f)_p^p\right)^{1/p}\leq \theta e_{T}(f),
$$
we bisect $T$ towards the edge $e$ (greedy bisection).  Otherwise, we bisect $T$ 
from its most recently generated vertex (newest vertex bisection).
Here $\theta$ is a fixed number in $(0,1)$.
In the case $p=\infty$ we use the condition
$$
\max\{e_{T_e^1}(f),e_{T_e^2}(f)\} \leq \theta e_{T}(f).
$$
This new refinement rule benefits from the mesh size reduction properties
of newest vertex bisection. Indeed, a property illustrated in Figure 1 
is that a sequence $\{BNN\}$ of one arbitrary bisection ($B$) followed 
by two newest vertex bisections ($N$) produces triangles with diameters
bounded by half the diameter of the initial triangle. A more general
property - which proof is elementary yet tedious - is the following: 
a sequence of the type $\{BNB \cdots BN\}$ of 
length $k+2$ with a newest vertex bisection at iteration $2$ and $k+2$ 
produces triangles with diameter 
bounded by  $(1-2^{-k})$ times the diameter of the initial triangle,
the worst case being illustrated in Figure 2 with $k=3$.

\begin{figure}[htbp]
\centerline{
\psfig{figure=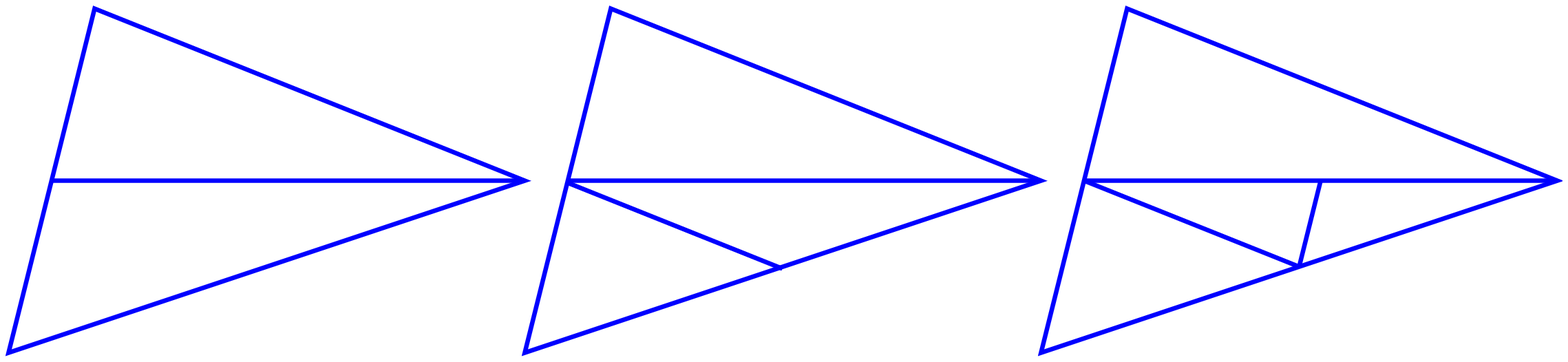,width=9cm}
}
\centerline{Figure 1: diameter reduction by a $\{BNN\}$ sequence}
\end{figure}

\begin{figure}[htbp]
\centerline{
\psfig{figure=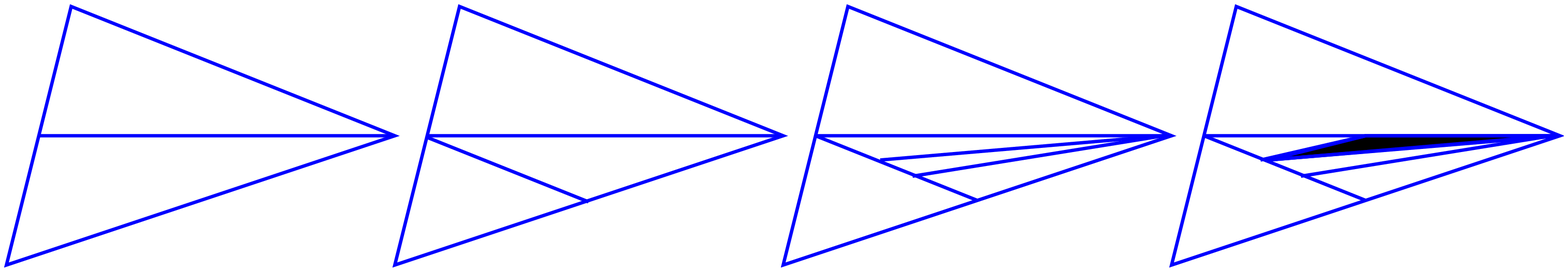,width=12cm}
}
\centerline{Figure 2: diameter reduction by a $\{BNBBN\}$ sequence}
\centerline{(the dark
triangle has diameter at most $7/8$ times the initial diameter)}
\end{figure}

Our next result shows that the modified algorithm now converges
for any $f\in L^p$.

\begin{theorem}
With $\cR$ defined as the modified bisection rule, we have
$$
f\in V(f,\cR)_p,
$$
for any $f\in L^p(\Omega)$ (or $\cC(\Omega)$ when $p=\infty$).
\end{theorem}
\noindent
{\bf Proof:} We first give the proof when $p<\infty$.
For each triangle $T\in \cD_j$ with $j\geq 1$, we introduce 
the two quantities 
$$
\alpha(T):=\frac {e_T(f)^p_p}{e_T(f)^p_p+e_{T'}(f)^p_p}\;\;{\rm and}
\;\; \beta(T):=\frac {e_T(f)^p_p+e_{T'}(f)^p_p}{e_{\cP(T)}(f)^p_p},
$$
where $T'$ is the ``brother'' of $T$, i.e. $\cC(\cP(T))=\{T,T'\}$.
When a greedy bisection occurs in the split of $\cP(T)$, we have
$\beta(T)\leq \theta^p$. When a newest vertex bisection occurs, we have
$$
\beta(T)\leq  C^p \frac {\inf_{\pi\in\Pi_m}\|f-\pi\|_{L^p(T)}^p+\inf_{\pi\in\Pi_m}\|f-\pi\|_{L^p(T')}^p}
{\inf_{\pi\in\Pi_m}\|f-\pi\|_{L^p(\cP(T))}^p}\leq C^p,
$$
where $C$ is the constant of \iref{unileb}.

We now consider a given level index $j>0$ 
and the triangulation $\cD_j$. For each $T\in \cD_j$,
we consider the chain of nested triangles
$(T_n)_{n=0}^j$ with $T_j=T$ and $T_{n-1}=\cP(T_n)$,
$n=j,j-1,\cdots,1$. 
We define
$$
\o\alpha(T)=\prod_{n=1}^j\alpha(T_n)
\quad \text{ and } \quad
\o\beta(T)=\prod_{n=1}^j\beta(T_n).
$$
It is easy to see that 
$$
\o\alpha(T)\o\beta(T)=\frac {e_T(f)^p_p}{e_{T_0}(f)^p_p}
$$
so that 
$$
e_T(f)^p_p\leq C_0\o\alpha(T)\o\beta(T),
$$
with $C_0:=\max_{T_0\in\cD_0}e_{T_0}(f)^p_p$.
It is also easy to check by induction on $j$ that
$$
\sum_{T\in\cD_j}\o\alpha(T)=\sum_{T\in\cD_{j-1}}\o\alpha(T)=\cdots=\sum_{T\in\cD_{1}}\o\alpha(T)=\#\cD_0.
$$
We denote by $f_j$ the approximation to $f$ in $V_j$ defined
by $f_j=\cA_T f$ on all $T\in\cT_j$ so that
$$
\|f-f_j\|_{L^p}^p=\sum_{T\in\cD_j}e_T(f)^p_p.
$$
In order to prove that $f_j$ converges to $f$
in $L^p$, it is sufficient to show that the sequence
$$
\e_j:=\max_{T\in\cD_j}\min\{\o\beta(T),{\rm diam}(T)\},
$$
tends to $0$ as $j$ grows. Indeed, if this holds, we
split $\cD_j$ into two sets $\cD_{j}^+$ and $\cD_j^-$ over
which $\o\beta(T)\leq \e_j$ and ${\rm diam}(T)\leq \e_j$
respectively. We can then write
\begin{eqnarray*}
\|f-f_j\|_{L^p}^p  & = & \sum_{T\in\cD_j^+}e_T(f)^p_p+\sum_{T\in\cD_j^-}e_T(f)^p_p \\
& \leq & C_0\e_j\sum_{T\in\cD_j^+}\o\alpha(T)+C^p\sum_{T\in\cD_j^-}\inf_{\pi\in\Pi_m}\|f-\pi\|_{L^p(T)}^p \\
& \leq & C_0\#\cD_0\e_j+C^p\sum_{T\in\cD_j^-}\inf_{\pi\in\Pi_m}\|f-\pi\|_{L^p(T)}^p,
\end{eqnarray*}
where $C$ is the constant of \iref{unileb}. Clearly the first term tends to $0$ and
so does the second term by standard properties of $L^p$
spaces since the diameter of the triangles in $\cD_j^-$ goes to $0$.
It thus remains to prove that
$$
\lim_{j\to +\infty}\e_j=0.
$$
Again we consider the chain $(T_n)_{n=0}^j$ which terminates at $T$,
and we associate to it a chain $(q_n)_{n=0}^{j-1}$ where
$q_n=1$ or $2$ if bisection of $T_n$ is greedy or newest vertex respectively.
If $r$ is the total number of $2$ in the chain $(q_n)$ we have
$$
\o\beta(T)\leq C^{pr}\theta^{p(j-r)},
$$
with $C$ the constant in \iref{unileb}.
Let a $k>0$ be a fixed number, large enough such that
$$
C \theta^{k-1} \leq 1.
$$
We thus have
\be
\o\beta(T)\leq (C^p \theta^{p(k-1)})^r\theta^{p(j-rk)}\leq \theta^{p(j-rk)}.
\label{controlbeta}
\ee
We now denote by $l$ the maximal number of disjoint 
sub-chains of the type $(\nu_1,2,\nu_2,\nu_3,\cdots,\nu_{q},2)$ 
with $\nu_j\in\{1,2\}$ and of length $q+2 \leq 2k+3$ which can be extracted
from $(q_n)_{n=0}^{j-1}$. From the remarks on the diameter reduction
properties of newest vertex bisection, we see that
$$
{\rm diam}(T) \leq B(1-2^{-2k})^l,
$$
with $B:=\max_{T_0\in\cD_0}{\rm diam}(T_0)$ a fixed constant. On the other
hand, it is not difficult to check that
\be
r\leq 3l+3+\frac{j-r}{2k}.
\label{lrjk}
\ee
Indeed let $\alpha_0$ be the total number of $1$
in the sequence $(q_n)$ which are not preceeded by a $2$,
and let $\alpha_i$ be the size of the series
of $1$ following the $i$-th occurence of $2$ in $(q_n)$ for
$i=1,\cdots,r$. Note that some $\alpha_i$ might be $0$.
Clearly we have
$$
j = \alpha_0+ \alpha_1+\cdots+\alpha_r+r.
$$
From the above equality, the number of $i$ such that $\alpha_i>2k$ is 
less than $\frac {j-r}{2k}$ and therefore there
is at least $m\geq r-\frac {j-r}{2k}$ indices
$\{i_0,\cdots,i_{m-1}\}$ such that $\alpha_i\leq 2k$.
Denoting by $\beta_i$ the position of the
$i$-th occurence of $2$ (so that $\beta_{i+1}=\beta_i+\alpha_i+1$),
we now consider the 
disjoint sequences of indices
$$
S_t=\{\beta_{i_{3t}},\cdots,\beta_{i_{3t+2}}\}, \;\; t=0,1,\cdots
$$
There is at least $\frac m 3-1$ such sequences within $\{1,\cdots,j\}$
and by construction each of them contains 
a sequence of the type $(\nu_1,2,\nu_2,\nu_3,\cdots,\nu_{q},2)$ 
with $\nu_j\in\{1,2\}$ and of length $q+2 \leq 2k+3$. Therefore
the maximal number of disjoint sequences of such type satisfies
$$
l\geq \frac 1 3(r-\frac {j-r}{2k})-1,
$$
which is equivalent to \iref{lrjk}. Therefore according to \iref{controlbeta}
$$
\o\beta(T)\leq \theta^{p(j-rk)} \leq \theta^{p(\frac j 2-3(l+1)k+\frac r 2)}
$$
If $3(l+1)k\leq \frac j 4$, we have 
$$
\o\beta(T)\leq \theta^{\frac {pj} 4}.
$$
On the other hand, if $3(l+1)k\geq \frac j 4$, we have
$$
{\rm diam}(T) \leq B(1-2^{-2k})^{\frac j{12 k}-1}.
$$
We therefore conclude that $\e_j$ goes to $0$ as $j$ grows, which proves the result
for $p<\infty$. 
\nl
\nl
We briefly sketch the proof for $p=\infty$, which is simpler. We now define $\beta(T)$ as
$$
\beta(T):=\frac {e_T(f)_\infty}{e_{\cP(T)}(f)_\infty},
$$
so that $\beta(T)\leq \theta$ if a greedy bisection occurs in the split of $\cP(T)$.
With the same definition of $\o\beta(T)$ we now have
$$
e_T(f)_\infty \leq C_0\o\beta(T),
$$
where $C_0:=\max_{T_0\in\cD_0} e_{T_0}(f)_\infty$. With the same definition of 
$\e_j$ and splitting
of $\cD_j$, we now reach
$$
\|f-f_j\|_{L^\infty}\leq \max\left \{C_0\e_j, C\max_{T\in\cD_j^-}\inf_{\pi\in\Pi_m}\|f-\pi\|_{L^\infty(T)}\right \},
$$
which again tends to $0$ if $\e_j$ tends to $0$ and $f$ is continuous.
The proof that $\e_j$ tends to $0$ as $j$ grows is then similar to the
case $p<\infty$.
\hfill $\diamond$

\begin{remark}
The choice of the parameter $\theta<1$ deserves some attention: if it is chosen
too small, then most bisections are of type $N$ and we end up with an isotropic
triangulation. In the case $m=1$,
a proper choice can be found by observing that when 
$f$ has $C^2$ smoothness, it can be locally approximated
by a quadratic
polynomial $q\in \Pi_2$ with $e_T(f)_p\approx e_T(q)_p$
when the triangle $T$ is small enough.
For such quadratic functions $q\in\Pi_2$, one can
explicitely study the minimal error reduction which is always ensured by the
greedy refinement rule defined a given decision function.
In the particular case $p=2$ and with the choice $\cA_T=P_T$
which is considered in the numerical experiments of \S 4,
explicit formulas for the error $\|q-P_Tq\|_{L^2(T)}$ 
can be obtained by formal computing and can be used to prove a guaranteed
error reduction by a factor $\theta^*=\frac 3 5$.
It is therefore natural to choose $\theta$ such that $\theta^*<\theta <1$
(for example $\theta=\frac 2 3$)
which ensures that bisections of type $N$ only occur in the early steps of the
algorithm, when the function still exhibits too many oscillations on 
some triangles.
\end{remark}

\section{Numerical illustrations}

The following numerical experiments were conducted with
piecewise linear approximation in the $L^2$ setting:
we use the $L^2$-based decision function
$$
d_T(f,e):=\|f-P_{T_e^1}f\|_{L^2(T_e^1)}^2+\|f-P_{T_e^2}f\|_{L^2(T_e^2)}^2.
$$
and we take for $\cA_T$ the $L^2(T)$-orthogonal projection $P_T$
onto $\Pi_1$. In these experiments, the function $f$
is either a quadratic polynomial or a function with a simple analytic expression
which allows us to compute the quantities $e_T(f)_2$ and $d_{T}(e,f)$
without any quadrature error,  or a numerical image in which case
the computation of these quantities is discretized on the pixel grid.

\subsection{Quadratic functions}

Our first goal is to illustrate numerically the optimal
adaptation properties of the refinement procedure 
in terms of triangle shape.
For this purpose, we take $f=\bq$ a quadratic form
i.e. an homogeneous polynomial of degree $2$. In this case,
all triangles should have the same aspect ratio 
since the Hessian is constant. In order to measure the
quality of the shape of a triangle $T$ in relation to $\bq$, we introduce
the following quantity: if $(a,b,c)$ are the edge vectors of $T$, we define
$$
\rho_\bq(T) := \frac{\max \{|\bq(a)|,|\bq(b)|,|\bq(c)|\}}{|T| \sqrt{|{\rm det}(\bq)|}},
$$
where ${\rm det}(\bq)$ is the determinant of the $2\times 2$ symmetric matrix $Q$
associated with $\bq$, i.e. such that 
$$
\bq(u)=\<Qu,u\>
$$
for all $u\in\R^2$. Using the reference triangle
and an affine change of variables, it is proved in \S2 of \cite{CM} that 
$$
e_T(\bq)_p\sim |T|^{1+\frac 1 p}\rho_\bq(T)\sqrt{|{\rm det}(\bq)|},
$$
with equivalence constants independent of $\bq$ and $T$.
Therefore, if $T$ is a triangle of given area, its shape should be
designed in order to minimize $\rho_\bq(T)$.

In the case where $\bq$ is positive definite or negative
definite, $\rho_\bq(T)$ takes small values when $T$ is isotropic
{\it with respect to the metric} $|(x,y)|_{\bq}:=\sqrt{|\bq(x,y)|}$, the minimal value $\frac 4 {\sqrt{3}}$ being
attained for an equilateral triangle for this metric.
Specifically, we choose $\bq(x,y):=x^2+100 y^2$ and
display in Figure 3 (left) the triangulation $\cD_{8}$ obtained after $j=8$ iterations
of the refinement procedure, starting with a triangle which is equilateral for 
the euclidean metric (and therefore not adapted to $\bq$). 
Triangles such that $\rho_\bq(T)\leq 4\sqrt 3$ (at most $3$ times the minimal value)
are displayed in white, 
others in grey. We observe that most triangles produced by
the refinement procedure are of the first type and therefore
have a good aspect ratio.

\begin{figure}[htbp]
\centerline{
\psfig{figure=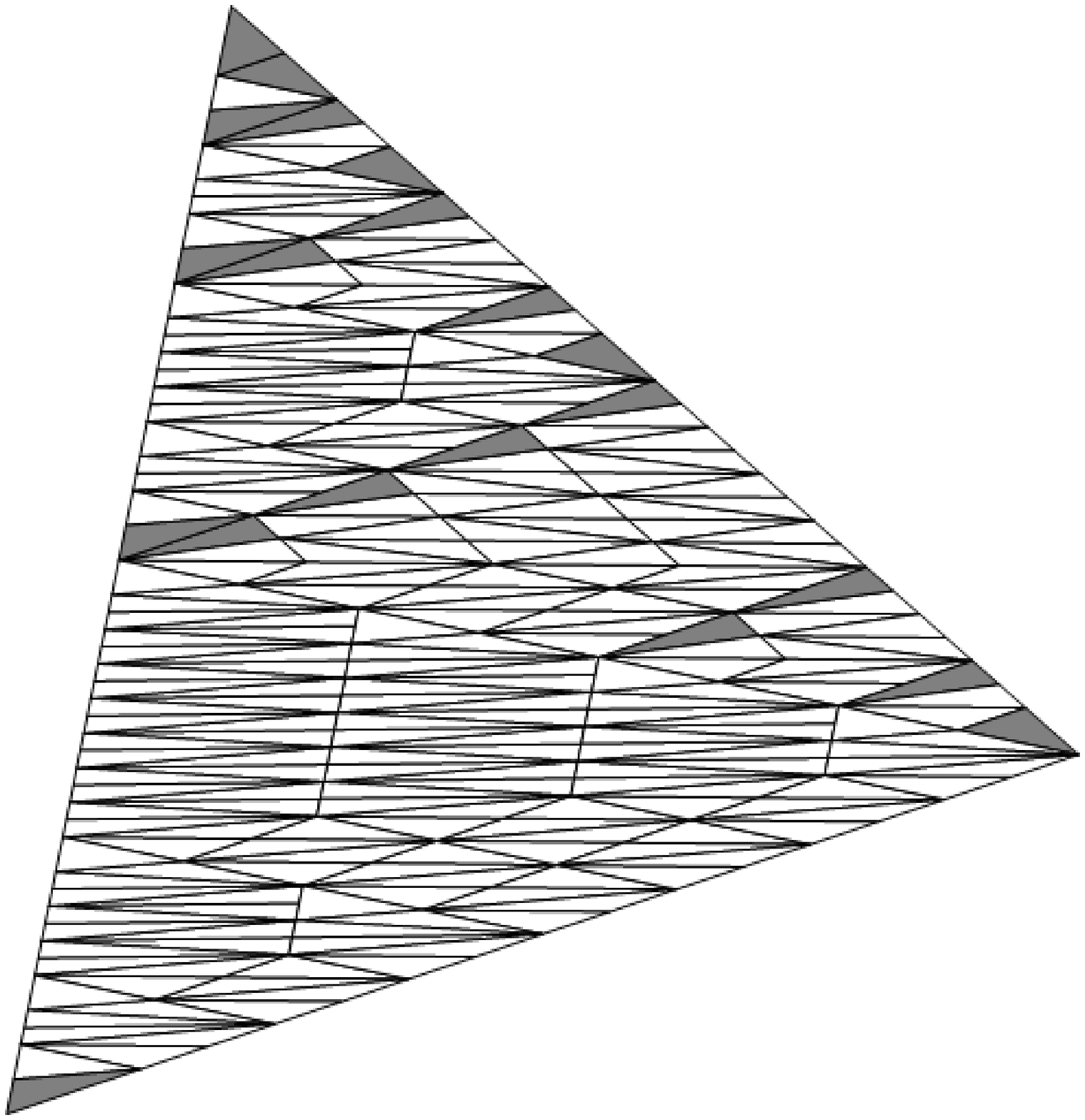,width=4cm,height=4cm}
\hspace{1cm}
\psfig{figure=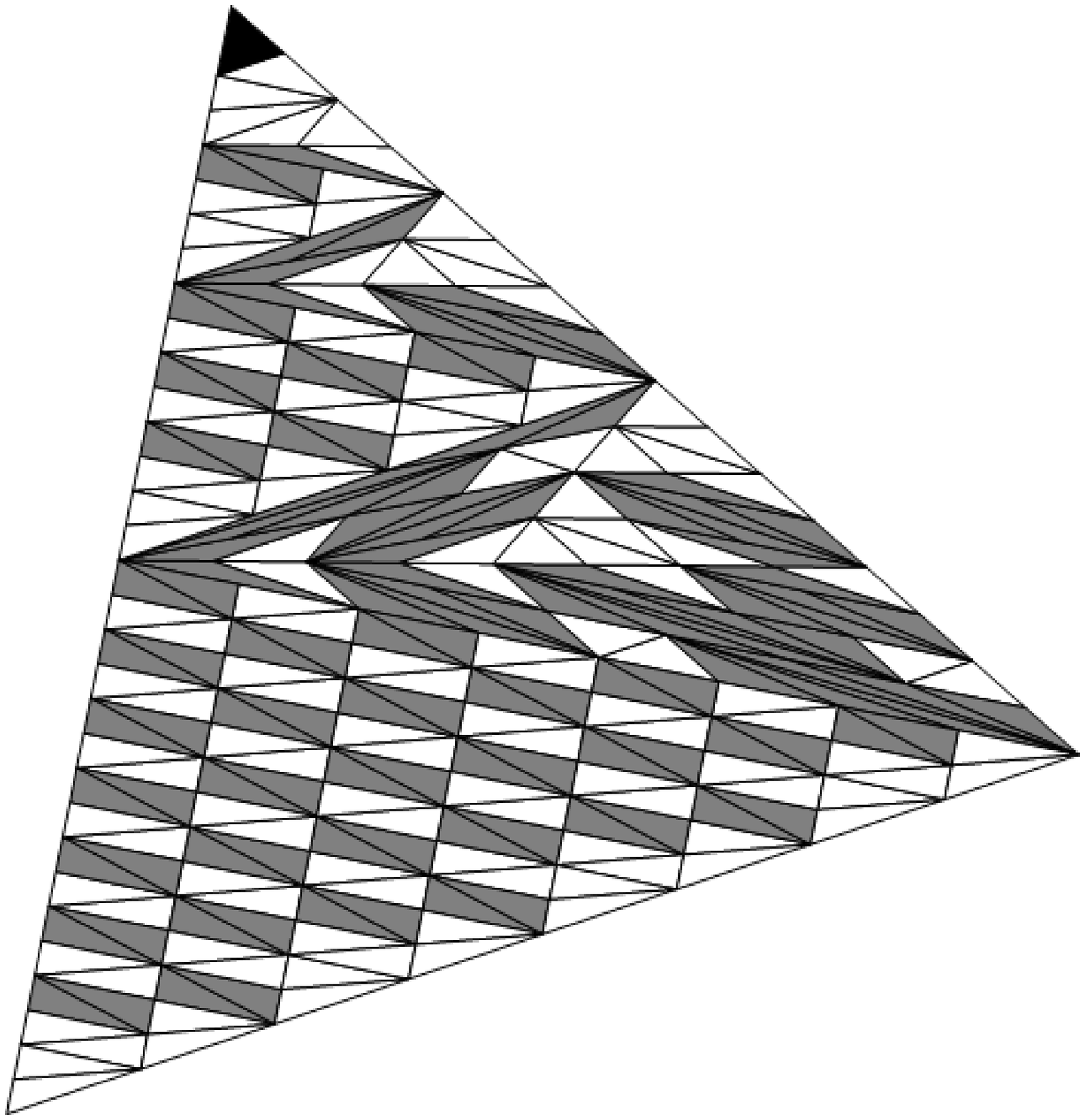,width=4cm,height=4cm}
}
\centerline{Figure 3: $\cD_{8}$ for $q(x,y):=x^2+100 y^2$ (left) and $q(x,y):=x^2-10 y^2$ (right).}
\end{figure}

The case of a
quadratic function of mixed signature is illustrated 
in Figure 3 (right) with $\bq(x,y):=x^2-10 y^2$. 
For such quadratic functions, triangles which are isotropic with
respect to the metric $|\cdot|_{|\bq|}$ have a low value of $\rho_\bq$,
where $|\bq|$ denotes the positive quadratic form associated
to the absolute value $|Q|$ of the symmetric matrix $Q$ associated
to $\bq$.
Recall for any symmetric matrix $Q$ there exists $\lambda_1,\lambda_2\in \R$ and a rotation $R$ such that 
$$
Q = R^\trans 
\left(\begin{array}{cc}
\lambda_1 & 0\\
0 & \lambda_2
\end{array}\right)
R,
$$
and the absolute value $|Q|$ is defined as 
$$
|Q| = R^\trans
\left(\begin{array}{cc}
|\lambda_1| & 0\\
0 & |\lambda_2|
\end{array}\right)
R.
$$
In the present case, $R=I$ and $|\bq|(x,y)=x^2+10 y^2$.

But one can also check that 
$\rho_\bq$ is left invariant by any linear transformation
with eigenvalues $(t,\frac 1 t)$ for any $t>0$ and eigenvectors $(u,v)$ such that
$\bq(u)=\bq(v)=0$, i.e. belonging to 
the \emph{null cone} of $\b q$. More precisely, for any such transformation
$\psi$ and any triangle $T$, one has $\rho_{\bq}(\psi(T))=\rho_\bq(T)$.
In our example we have $u=(\sqrt {10},1)$ and $v=(\sqrt{10},-1)$).
Therefore long and thin triangles which are aligned with
these vectors also have a low value of $\rho_\bq$.
Triangles $T$
such that $\rho_{|\bq|}(T)\leq 4 \sqrt 3$ are displayed in white, those 
such that $\rho_{\bq}(T)\leq 4\sqrt 3$ while $\rho_{|\bq|}(T)> 4 \sqrt 3$ -
i.e. adapted to $\bq$ but not to $|\bq|$ - are displayed in grey,
and the others in dark. We observe that all the triangles triangles produced by
the refinement procedure except one are either of the first or second type
and therefore have a good aspect ratio.

\subsection{Sharp transition}

We next study the adaptive triangulations produced by the greedy tree
algorithm for a function $f$ displaying a sharp transition along a curved edge. Specifically
we take
$$
f(x,y)=f_\delta(x,y) := g_\delta(\sqrt{x^2+y^2}),
$$ 
where $g_\delta$
 is defined by
$g_\delta(r)= \frac{5-r^2} 4$ for $0\leq r\leq 1$, $g_\delta(1+\delta+r)=-\frac{5-(1-r)^2} 4$ for $r\geq 0$, and $g_{\delta}$ is a polynomial of degree $5$ on $[1,1+\delta]$
which is determined by imposing that
$g_\delta$ is globally $C^2$.
The parameter $\delta$ therefore measures the sharpness of the transition as illustrated
in Figure 4.
It can be shown that the Hessian of $f_\delta$ is negative definite for $\sqrt{x^2+y^2}< 1+\delta/2$, and
of mixed type for $1+\delta/2 < \sqrt{x^2+y^2}\leq 2+\delta$.

\begin{figure}[htbp]
\centerline{
\psfig{figure=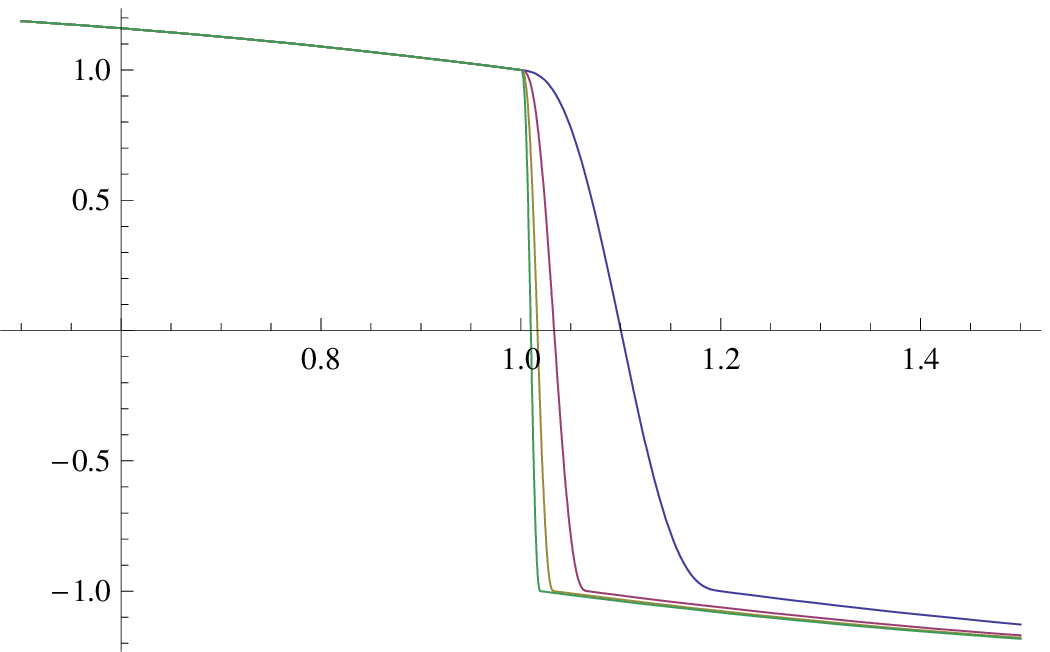,width=4cm,height=3cm}
}
\centerline{Figure 4: The function $g_\delta$, for $\delta=0.02,0.03,0.07,0.2$.}
\end{figure}

Figure 5 displays the triangulation $\cT_{10000}$ obtained after $10000$ steps of
the algorithm for $\delta=0.2$. In particular, triangles $T$ such that $\rho_\bq(T)\leq 4$ 
- where $\bq$ is the quadratic form associated with $d^2f$ measured at the barycenter of $T$ - 
are displayed in white, others in grey. As expected, most triangles are of the
first type and therefore well adapted to $f$. We also display
on this figure the adaptive isotropic triangulation produced by the greedy tree algorithm
based on newest vertex bisection for the same number of triangles.

\begin{figure}[htbp]
\centerline{
\psfig{figure=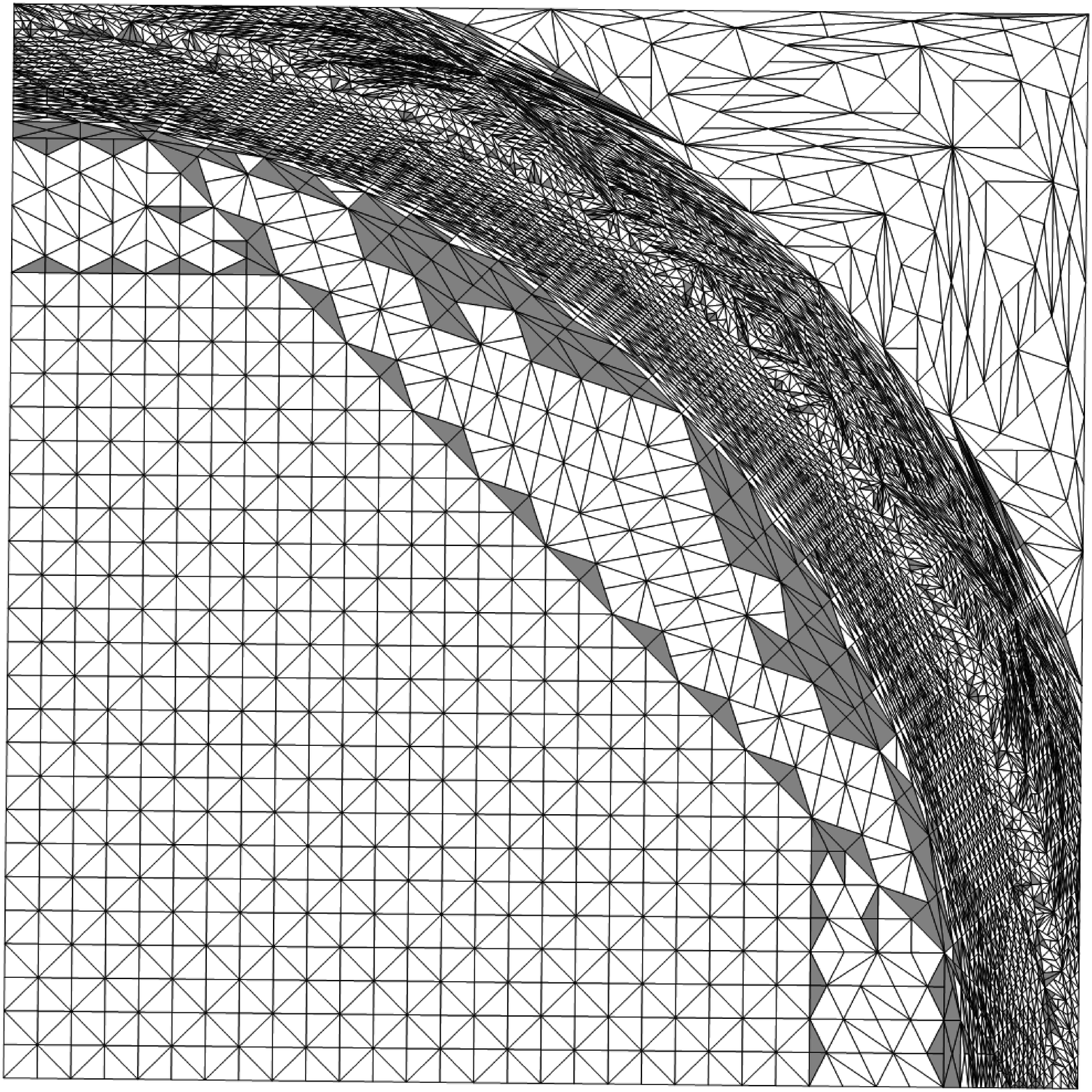,width=4cm,height=4cm}
\psfig{figure=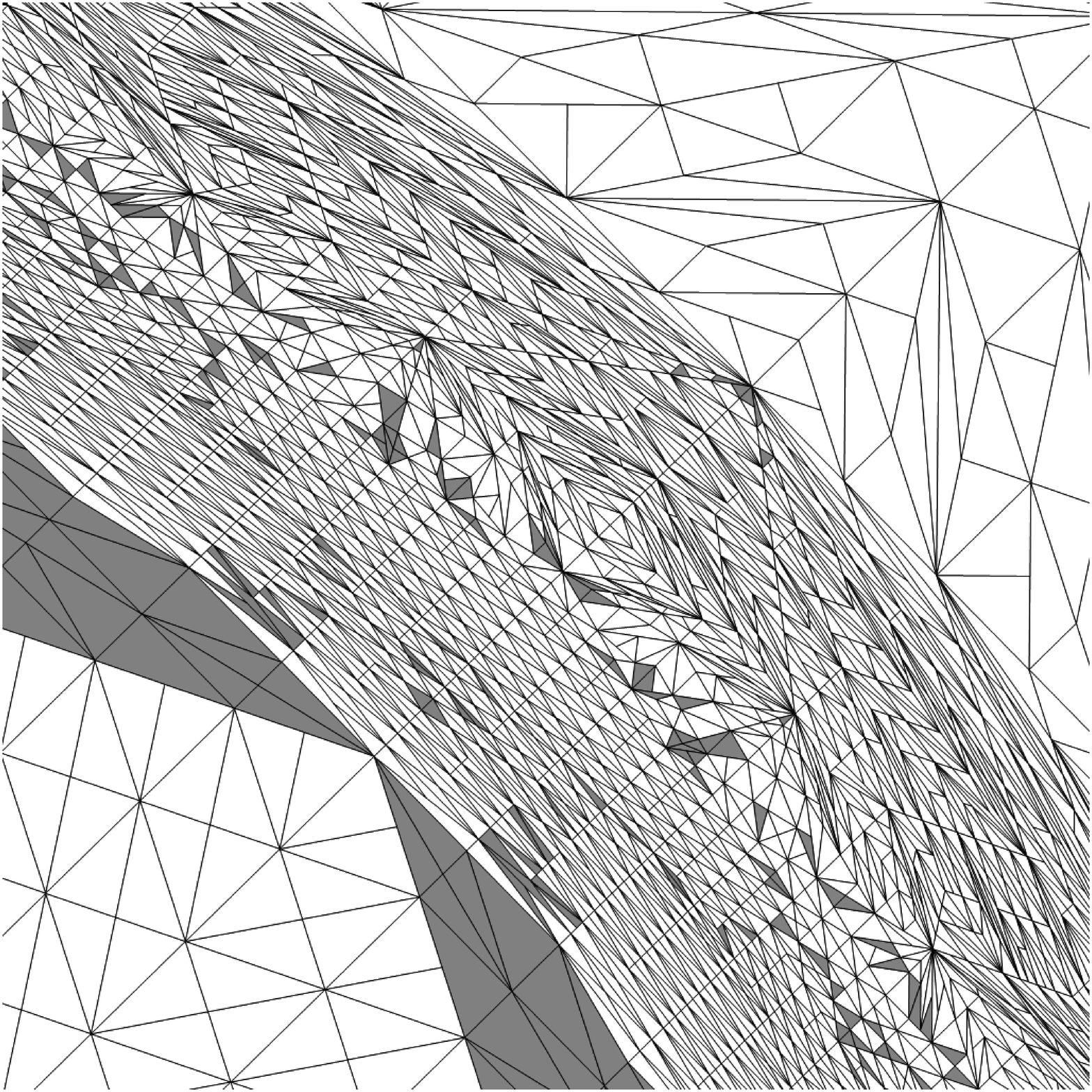,width=3.9cm,height=3.9cm}
\psfig{figure=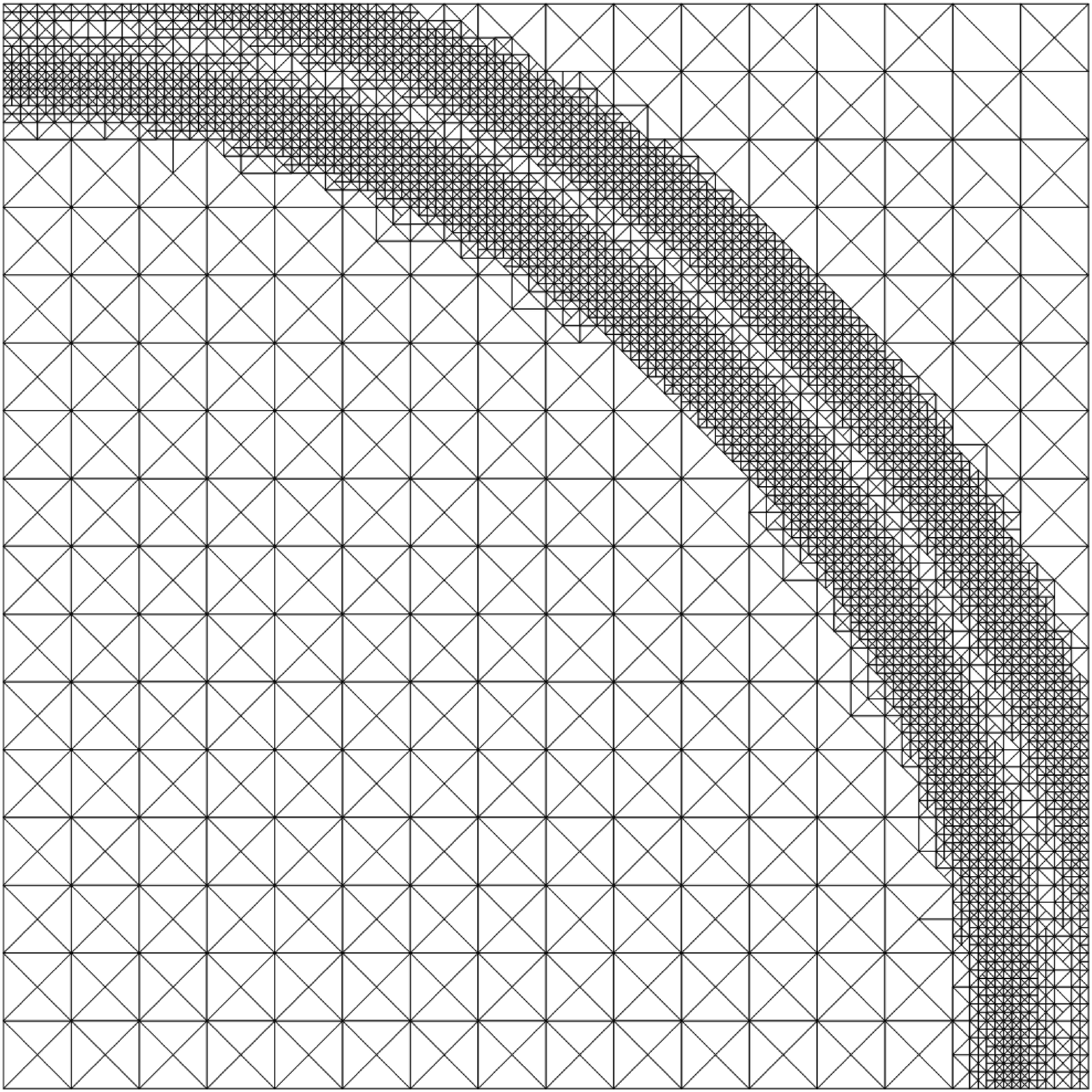,width=4cm,height=4cm}}
\centerline{Figure 5:   $\cT_{10000}$ (left),  detail (center), isotropic triangulation (right).}
\end{figure}

Since $f$ is a $C^2$ function, approximations by uniform, adaptive isotropic and adaptive anisotropic triangulations all yield the convergence rate $\cO(N^{-1})$.
However the constant 
$$
C:=\limsup_{N\to +\infty} N\|f-f_N\|_{L^2},
$$ 
strongly differs depending on the algorithm
and on the sharpness of the transition, as illustrated in the table below.
We denote by $C_U$, $C_I$ and $C_A$ the empirical constants
(estimated by $N \|f-f_N\|_2$ for $N=8192$)
in the uniform, adaptive isotropic and adaptive anisotropic case respectively, and by
$U(f):=\|d^2f\|_{L^{2}}$,
$I(f):=\|d^2f\|_{L^{2/3}}$ and
$A(f):=\|\sqrt{|{\rm det}(d^2f)|}\|_{L^{2/3}}$
the theoretical constants suggested by
\iref{unier}, \iref{isoer}
and \iref{aniser}. We observe that $C_U$
and $C_I$ grow in a similar way as $U(f)$ and $I(f)$
as $\delta\to 0$ (a detailed computation shows that $U(f)\approx 10.37 \, \delta^{-3/2}$
and $I(f)\approx 14.01 \, \delta^{-1/2}$).
In contrast $C_A$ and $A(f)$ remain uniformly bounded, a fact
which reflects the superiority of the anisotropic mesh
as the layer becomes thinner.

$$
\begin{array}{c|c|c|c|c|c|c|}
\delta & U(f)  & I(f) & A(f) & C_U & C_I   & C_A
\\
\hline
0.2	& 103  & 27 & 6.75 & 7.87	& 1.78 &  0.74\\
0.1 & 602 	& 60 & 8.50 & 23.7	& 2.98 & 0.92\\
0.05& 1705 & 82 & 8.48 & 65.5	&	4.13 & 0.92\\
0.02& 3670 & 105 & 8.47 & 200	&	6.60 & 0.92
\end{array}
$$

\subsection{Numerical images}

We finally apply the greedy tree
algorithm to numerical images. In this case
the data $f$ has the form of a discrete array of pixels,
and the $L^2(T)$-orthogonal projection is replaced
by the $\ell^2(S_T)$-orthogonal projection,
where $S_T$ is the set of pixels with centers contained in $T$.
The approximated $512\times 512$ image is displayed
in Figure 6 which also shows its approximation $f_{N}$
by the greedy tree algorithm
based on newest vertex bisection with $N=2000$ triangles.
The systematic use of isotropic triangles results in
strong ringing artifacts near the edges.

\begin{figure}[htbp]
\centerline{
\psfig{figure=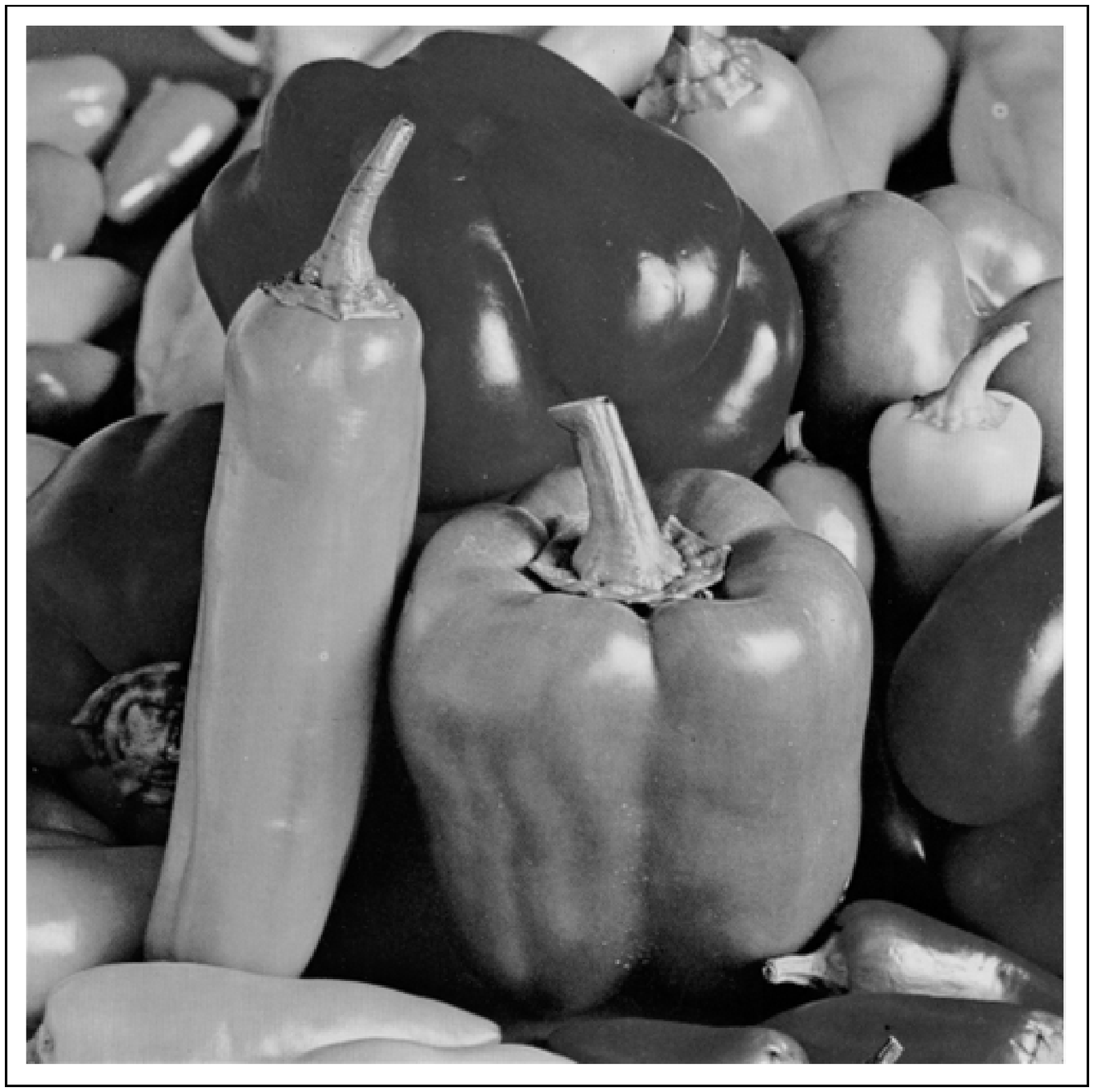,width=6cm,height=6cm}
\psfig{figure=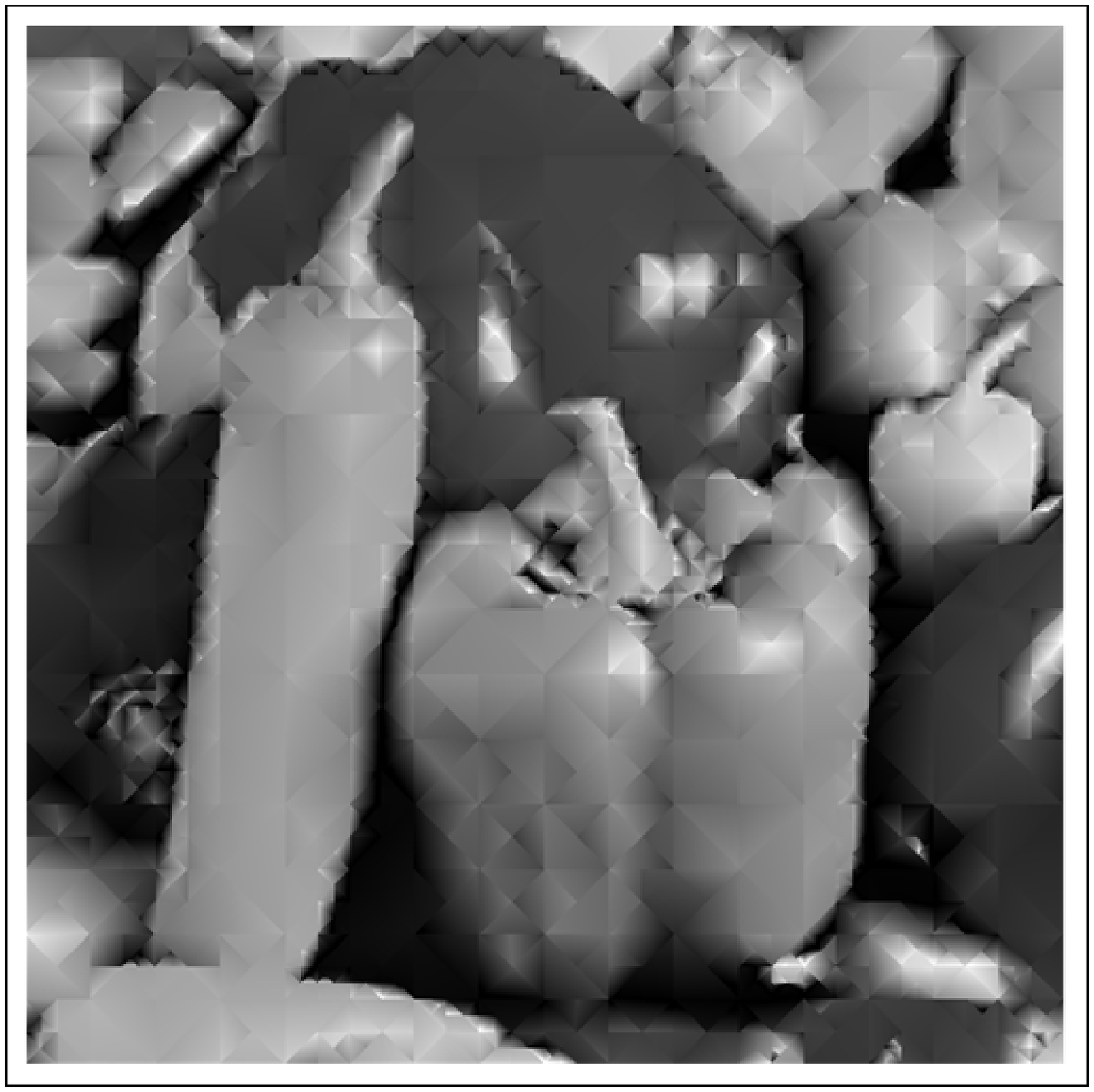,width=6cm,height=6cm}
}
\centerline{Figure 6:   The image ''peppers'' (left),  $f_{2000}$ with newest vertex (right).}
\end{figure}

\begin{figure}[htbp]
\centerline{
\psfig{figure=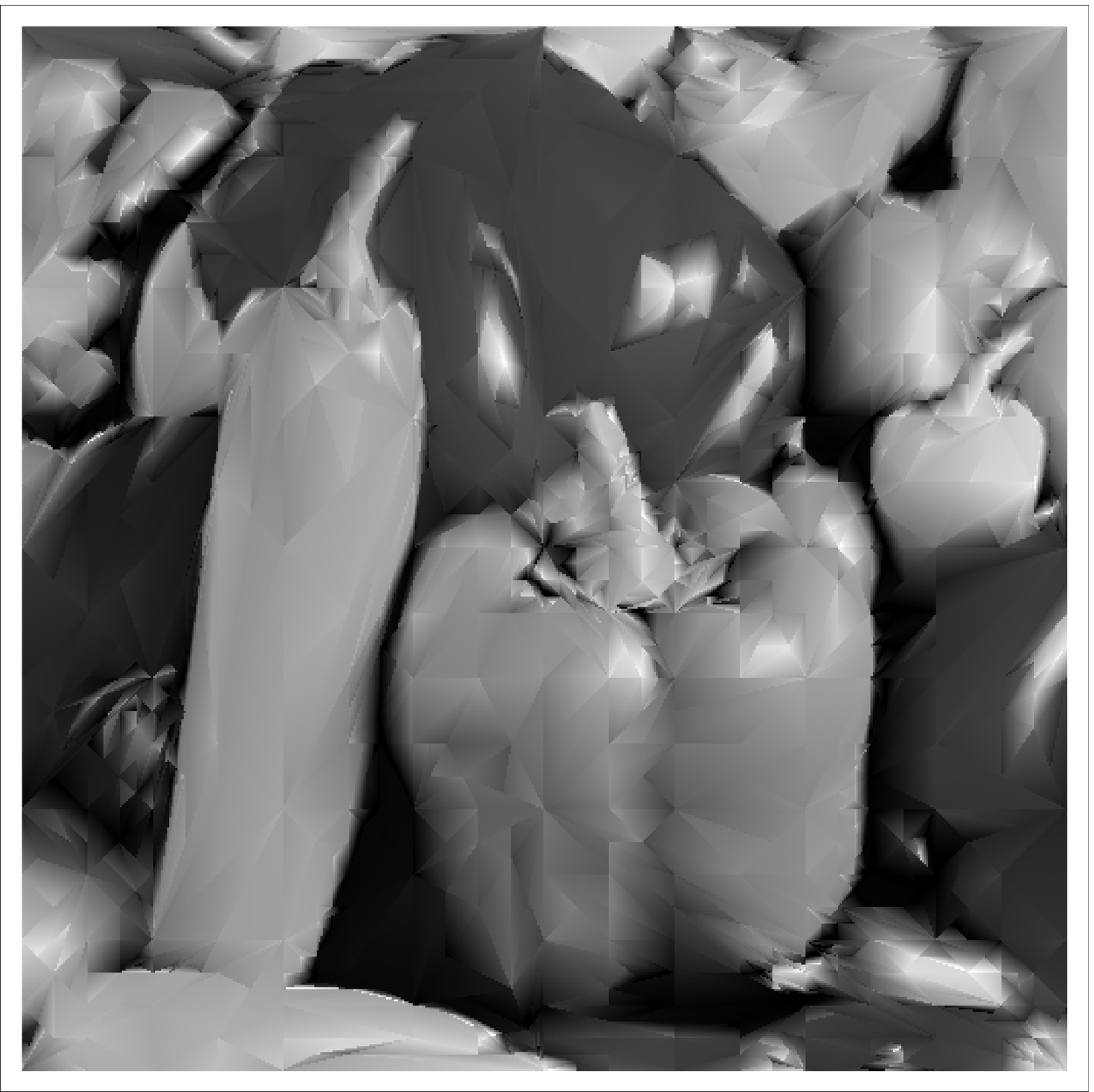,width=6cm,height=6cm}
\psfig{figure=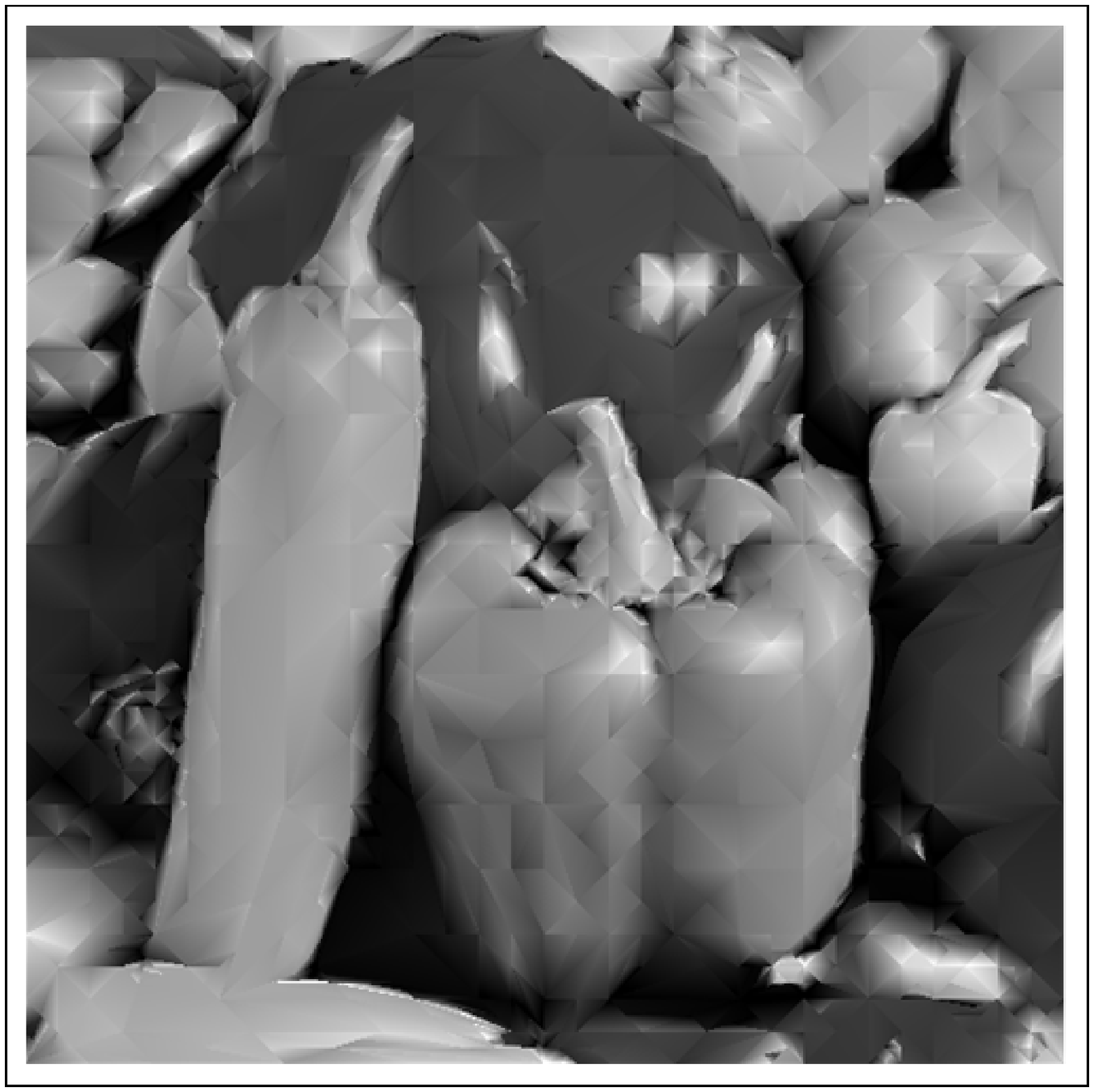,width=6cm,height=6cm}
}
\centerline{Figure 7:  $f_{2000}$ with greedy bisection (left),  modified procedure (right).}
\end{figure}

We display in Figure 7 the result of the same algorithm 
now based on our greedy bisection procedure 
with the same number of triangles. As expected, the edges
are better approximated due to the presence of well oriented
anisotropic triangles. Yet artifacts persist on certain edges
due to oscillatory
features in the image which tend to mislead the algorithm
in its search for triangles with good aspect ratio, as explained
in \S 3.2. These artifacts tend to disappear if we use the 
modified refinement rule proposed in \S 3.3 as also illustrated on
Figure 7. This modification is thus useful in the
practical application of the algorithm, in addition of being
necessary for proving convergence convergence of the
approximations towards any $L^p$ function. Note that
encoding a triangulation resulting from $N$ iterations 
of the anisotropic refinement algorithm is more costly than for the newest
vertex rule: the algorithm encounters at most $2N$ triangles
and for each of them, one needs to encode one out of four options
(bisect towards edge $a$ or $b$ or $c$ or not bisect), therefore
resulting into $4N$ bits,
while only two options need to be encoded when using the newest vertex rule (bisect or not),
therefore resulting into $2N$ bits. In the perspective of
applications to image compression, another issue is the 
quantization and encoding of the piecewise affine function
as well as the treatment of the triangular visual artifacts that are inherent
to the use of discontinuous piecewise polynomials on 
triangulated domains. These issues will be discussed in 
a further work specifically dealing with image applications.

\section{Conclusions and perspectives}

In this paper, we have studied a simple greedy 
refinement procedure which 
generates triangles that tend to 
have an {\it optimal aspect ratio}.
This fact is rigorously proved in \cite{CM},
together with the optimal 
convergence estimate \iref{aniser}
for the adaptive triangulations constructed by the
greedy tree algorithm in the case where 
the approximated function $f$ is
$C^2$ and convex. Our numerical results
illustrate these properties.

In the present paper we also show
that for a general $f\in L^p$ the refinement procedure 
can be misled by 
oscillations in $f$, and that
this drawback may be circumvented by a 
simple modification of the refinement procedure.
This modification appears to be useful 
in image processing applications, as shown
by our numerical results.

Let us finally mention several perspectives that are raised
from our work, and that are the object of current investigation:
\begin{enumerate}
\item
Conforming triangulations: our algorithm inherently generates
hanging nodes, which might not be desirable in certain applications,
such as numerical discretization of PDE's where anisotropic elements
are sometimes used \cite{Apel}.
When using the greedy tree algorithm, an obvious way of avoiding this phenomenon is to
bisect the chosen triangle together with an adjacent triangle in order
to preserve conformity. However, it is no more clear that
this strategy generates optimal triangulations. In fact, we observed
that many inappropriately oriented triangles can be generated by this
approach. An alternative strategy is to apply the non-conforming greedy tree algorithm
until a prescribed accuracy is met, followed by an additional
refinement procedure in order to remove hanging nodes.
\item
Discretization and encoding: 
our work is in part motivated by applications to image and terrain data
processing and compression. In such applications
the data to be approximated is usually given in discrete form
(pixels or point clouds) and the algorithm can be adapted
to such data, as shown in our numerical image examples. 
Key issues which need to be dealt with are then 
(i) the efficient encoding of the approximations and of the triangulations
using the tree structure in a similar spirit as in \cite{CDDD} and
(ii) the removal of the triangular visual artifacts due to discontinuous piecewise polynomial
approximation by an appropriate post-processing step.
\item
Adaptation to curved edges: 
one of the motivation for the use of anisotropic triangulations is the 
approximation of functions with jump discontinuities along an edge.
For simple functions, such as characteristic functions of domains with smooth boundaries,
the $L^p$-error rate with an optimally adapted triangulation of $N$ elements
is known to be $\cO(N^{-\frac 2 p})$. This rate reflects an $\cO(1)$ error concentrated
on a strip of area $\cO(N^{-2})$ separating the curved edge
from a polygonal line. Our first investigations in this direction indicate that the
greedy tree algorithm based on our refinement procedure 
cannot achieve this rate, due to the fact that bisection does not offer enough
geometrical adaptation. This is in contrast with other splitting procedures,
such as in \cite{DL} in which the direction of the new cutting edge is optimized within 
an infinite range of possible choices, or \cite{D} where the number of choices
grows together with the resolution level. An interesting question, addressed and partially answered in \S 9.1 of \cite{thesisJM},
is thus to understand if the optimal rate for edges can be achieved
by a splitting procedure with a small and fixed number 
of choices similar to our refinement procedure, which would be beneficial
from both a computational and encoding viewpoint.  
\end{enumerate}

$\;$
\nl
Albert Cohen
\nl
UPMC Univ Paris 06, UMR 7598, Laboratoire Jacques-Louis Lions, F-75005, Paris, France
\nl
CNRS, UMR 7598, Laboratoire Jacques-Louis Lions, F-75005, Paris, France
\nl
cohen@ann.jussieu.fr
\nl
\nl
Nira Dyn
\nl
School of Mathematics, Tel Aviv University, Ramat Aviv, Israel
\nl
niradyn@math.tau.ac.il
\nl\nl
Fr\'ed\'eric Hecht
\nl
UPMC Univ Paris 06, UMR 7598, Laboratoire Jacques-Louis Lions, F-75005, Paris, France
\nl
CNRS, UMR 7598, Laboratoire Jacques-Louis Lions, F-75005, Paris, France
\nl
hecht@ann.jussieu.fr
\nl
\nl
Jean-Marie Mirebeau
\nl
UPMC Univ Paris 06, UMR 7598, Laboratoire Jacques-Louis Lions, F-75005, Paris, France
\nl
CNRS, UMR 7598, Laboratoire Jacques-Louis Lions, F-75005, Paris, France
\nl
mirebeau@ann.jussieu.fr

\end{document}